\documentclass[a4paper,11pt]{article}
\usepackage{array}
\usepackage{theorem}
\usepackage{amsmath,amscd,amssymb} 
\usepackage{latexsym}
\usepackage{epsfig}

\theorembodyfont{\sl}

\newtheorem{lemma}{Lemma}[section]

\newtheorem{proposition}[lemma]{Proposition}
\newtheorem{theorem}[lemma]{Theorem}
\newtheorem{corollary}[lemma]{Corollary}

\newtheorem{introthm}{Theorem}

\newcommand{\CC}{\mathbb C}

\newcommand{\HH}{\mathbb H}

\newcommand{\NN}{\mathbb N}

\newcommand{\PP}{\mathbb P}
\newcommand{\QQ}{\mathbb Q}
\newcommand{\RR}{\mathbb R}

\newcommand{\UU}{\mathbb U}
\newcommand{\VV}{\mathbb V}
\newcommand{\WW}{\mathbb W}

\newcommand{\ZZ}{\mathbb Z}

\newcommand{\cD}{\mathcal D}

\newcommand{\cF}{\mathcal F}

\newcommand{\cH}{\mathcal H}

\newcommand{\cS}{\mathcal S}

\newcommand{\cV}{\mathcal V}

\newcommand{\To}{\longrightarrow}
\newcommand{\Mapsto}{\mapstochar\longrightarrow}

\newcommand{\tensor}{\otimes}

\renewcommand{\Tilde}{\widetilde}
\renewcommand{\Hat}{\widehat}
\renewcommand{\Bar}{\overline}

\newcommand{\cross}{\times}

\newcommand{\imic}{\cong}

\newcommand{\GL}{\mathop{\mathrm {GL}}\nolimits}

\newcommand{\SL}{\mathop{\mathrm {SL}}\nolimits}
\newcommand{\SO}{\mathop{\mathrm {SO}}\nolimits}

\newcommand{\Orth}{\mathop{\null\mathrm {O}}\nolimits}

\newcommand{\Aut}{\mathop{\mathrm {Aut}}\nolimits}

\newcommand{\Hom}{\mathop{\mathrm {Hom}}\nolimits}

\newcommand{\Stab}{\mathop{\mathrm {Stab}}\nolimits}

\newcommand{\Lift}{\mathop{\mathrm {Lift}}\nolimits}
\newcommand{\Mat}{\mathop{\mathrm {Mat}}\nolimits}
\newcommand{\Lie}{\mathop{\mathrm {Lie}}\nolimits}
\newcommand{\codim}{\mathop{\mathrm {codim}}\nolimits}
\newcommand{\diag}{\mathop{\mathrm {diag}}\nolimits}

\renewcommand{\Im}{\mathop{\mathrm {Im}}\nolimits}
\newcommand{\orb}{\mathop{\mathrm {orb}}\nolimits}

\newcommand{\rank}{\mathop{\mathrm {rank}}\nolimits}

\newcommand{\cusp}{\mathop{\mathrm {cusp}}\nolimits}

\newcommand{\trace}{\mathop{\mathrm {tr}}\nolimits}

\newcommand{\id}{\mathop{\mathrm {id}}\nolimits}

\newcommand{\divv}{\mathop{\null\mathrm {div}}\nolimits}
\newcommand{\ord}{\mathop{\null\mathrm {ord}}\nolimits}
\newcommand{\bdT}{\mathbf T}

\newcommand{\bde}{\mathbf e}

\newcommand{\bdv}{\mathbf v}
\newcommand{\bdw}{\mathbf w}
\newcommand{\bdx}{\mathbf x}

\newcommand{\latt}[1]{{\langle{#1}\rangle}}

\newcommand{\Kthree}{\mathop{\mathrm {K3}}\nolimits}
\newcommand{\qedsymbol}{\mbox{$\Box$}}
\newcommand{\qed}{\unskip\nobreak\hfil\penalty50\hskip1em\hbox{}\nobreak
\hfill\qedsymbol\parfillskip=0pt\finalhyphendemerits=0}
\newenvironment{proof}{\begin{ProofwCaption}{Proof}}{\end{ProofwCaption}}
\newenvironment{ProofwCaption}[1]
 {\addvspace\theorempreskipamount \noindent{\it #1.}\rm}
 {\qed \par \addvspace\theorempostskipamount}

\begin{document}

\title{The Kodaira dimension of the moduli of K3 surfaces}
\author{V. Gritsenko, K.~Hulek and G.K. Sankaran}
\maketitle
\begin{abstract}
  The global Torelli theorem for projective $\Kthree$ surfaces was
  first proved by Piatetskii-Shapiro and Shafarevich 35 years ago,
  opening the way to treat moduli problems for $\Kthree$ surfaces.
  The moduli space of polarised $\Kthree$ surfaces of degree $2d$ is a
  quasi-projective variety of dimension $19$. For general $d$ very
  little has been known about the Kodaira dimension of these
  varieties.  In this paper we present an almost complete solution to
  this problem.  Our main result says that this moduli space is of
  general type for $d>61$ and for $d=46$, $50$, $54$, $58$, $60$.
\end{abstract}

\section{Introduction}

Moduli spaces of polarised $\Kthree$ surfaces can be identified with
the quotient of a classical hermitian domain of type $IV$ and
dimension $19$ by an arithmetic group.  The general set-up for the
problem is the following.  Let $L$ be an integral lattice with a
quadratic form of signature $(2,n)$ and let
\begin{equation}\label{DL}
  \cD_L=\{[\bdw] \in \PP(L\otimes \CC) \mid 
  (\bdw,\bdw)=0,\ (\bdw,\Bar\bdw)>0\}^+
\end{equation}
be the associated $n$-dimensional Hermitian domain (here $+$ denotes
one of its two connected components). We denote by $\Orth(L)^+$ the
index $2$ subgroup of the integral orthogonal group $\Orth(L)$
preserving $\cD_L$.  We are, in general, interested in the birational
type of the $n$-dimensional variety
\begin{equation}\label{FLG}
  \cF_L(\Gamma)=\Gamma\backslash \cD_L
\end{equation}
where $\Gamma$ is a subgroup of $\Orth^+(L)$ of finite index. Clearly,
the answer will depend strongly on the lattice $L$ and the chosen
subgroup $\Gamma$.

A compact complex surface $S$ is a $\Kthree$ surface if $S$ is simply
connected and there exists a holomorphic $2$-form $\omega_S \in
H(S,\Omega^2)$ without zeros. For example, a smooth quartic in
$\PP^3(\CC)$ is a $\Kthree$ surface and all quartics (modulo
projective equivalence) form a (unirational) space of dimension $19$. 

The second cohomology group $H^2(S,\ZZ)$ with the intersection pairing
is an even unimodular lattice of signature $(3,19)$, more precisely,
\begin{equation}\label{LK3}
  H^2(S,\ZZ)\cong L_{\Kthree}=3U\oplus 2E_8(-1)
\end{equation}
where $U$ is the hyperbolic plane and $E_8(-1)$ is the negative
definite even lattice associated to the root system $E_8$.  The
$2$-form $\omega_S$, considered as a point of $\PP(L_{\Kthree}\otimes
\CC)$, is the period of $S$.  By the Torelli theorem the period of a
$\Kthree$ surface determines its isomorphism class. The moduli space
of all $\Kthree$ surfaces is not Hausdorff.  Therefore it is better to
restrict to moduli spaces of polarised $\Kthree$ surfaces.  The moduli
of all algebraic $\Kthree$ surfaces are parametrised by a countable
union of $19$-dimensional irreducible algebraic varieties. To choose a
component we have to fix a polarisation.  A polarised $\Kthree$
surface of degree $2d$ is a pair $(S,H)$ consisting of a $\Kthree$
surface $S$ and a primitive pseudo-ample divisor $H$ on $S$ of degree
$H^2=2d>0$. If $h$ is the corresponding vector in the lattice
$L_{\Kthree}$ then its orthogonal complement
\begin{equation}\label{L2d}
  h^{\perp}_{L_{\Kthree}}\cong L_{2d}=2U\oplus 2E_8(-1) 
  \oplus \langle -2d \rangle
\end{equation}
is a lattice of signature $(2,19)$. 

The $2$-form $\omega_S$ determines a point of $\cD_{L_{2d}}$ modulo
the group
$$
\Tilde{\Orth}^+(L_{2d}) = \{ g\in \Orth^+(L_{\Kthree}) \mid g(h)=h \}. 
$$
By the global Torelli theorem (\cite{P-SS}) and the surjectivity of
the period map
\begin{equation}\label{F2d}
  \cF_{2d}=\Tilde{\Orth}^+(L_{2d})\setminus \cD_{L_{2d}}
\end{equation}
is the coarse moduli space of polarised $\Kthree$ surfaces of degree
$2d$.  By a result of Baily and Borel~\cite{BB}, $\cF_{2d}$ is a
quasi-projective variety.  One of the fundamental problems is to
determine its birational type. 

For $d=2$, $3$ and $4$ the polarised $\Kthree$ surfaces of degree $2d$
are complete intersections in $\PP^{d+1}(\CC)$ and the moduli spaces
$\cF_{2d}$ for such $d$ are classically known.  Mukai has extended
these results in his papers \cite{Mu1}, \cite{Mu2} and \cite{Mu3} to
$1\leq d\leq 10$ and $d=17$, $19$, showing that these moduli spaces
are also unirational. 

In the other direction there are two results of Kondo and of
Gritsenko.  Kondo \cite{Ko1} considered the moduli spaces $\cF_{2p^2}$
where $p$ is a prime number.  (The reason for this choice is that all
these spaces are covers of $\cF_{2}$.)  He proved that these spaces
are of general type for $p$ sufficiently large.  His result, however,
is not effective.  Gritsenko \cite{G} showed a result for level
structures: let $\Tilde{\Orth}^+(L_{2d})(q)$ be the principal
congruence subgroup of $\Tilde{\Orth}^+(L_{2d})$ of level $q$. Then
$\Tilde{\Orth}^+(L_{2d})(q)\setminus \cD_{L_{2d}}$ is of general type
for any $d$ if $q\ge 3$.  In this paper we determine the Kodaira
dimension of $\cF_{2d}$ without imposing any {\it a priori\/}
restriction on $d$.  \smallskip

\begin{introthm}\label{mainthm}
  The moduli space $\cF_{2d}$ of $\Kthree$ surfaces with a
  polarisation of degree $2d$ is of general type for any $d>61$ and
  for $d=46$, $50$, $54$, $57$, $58$ and $60$. 

  If $d\ge 40$ and $d\ne 41$, $44$, $45$ or $47$ then the Kodaira
  dimension of $\cF_{2d}$ is non-negative. 
\end{introthm}

The description of the moduli space $\cF_{2d}$ as a quotient of the
symmetric space $\cD_{L_{2d}}$ by a subgroup of the orthogonal group
leads us to study, more generally, quotients of the form
$\cF_L(\Gamma)=\Gamma\backslash \cD_L$.  One of the main tools in our
proof of the main theorem is the following general result (for a more
precise formulation see Theorem~\ref{main_sings_theorem}). 

\begin{introthm}\label{sings_theorem_intro}
  Let $L$ be a lattice of signature $(2,n)$ with $n\ge 9$, and let
  $\Gamma<\Orth^+(L)$ be a subgroup of finite index. Then there exists
  a toroidal compactification $\Bar\cF_L(\Gamma)$ of
  $\cF_L(\Gamma)=\Gamma\backslash\cD_L$ such that $\Bar\cF_L(\Gamma)$
  has canonical singularities. 
\end{introthm}

We hope that this result will also be important for other applications. 

The plan of the paper is as follows. In Section~\ref{orthogonal} we
give the basic definitions that we shall need and explain what the
obstructions are to showing that $\cF_L(\Gamma)$ is of general type. 
These obstructions may be called elliptic, cusp and reflective. The
elliptic obstructions come from singularities of $\cF_L(\Gamma)$ and
its compactifications. The cusp obstructions come from infinity,
{i.e.} from the fact that $\cF_L(\Gamma)$ is only quasi-projective. 
The reflective obstructions come from divisors fixed by $\Gamma$ in
its action on the symmetric space $\cD_L$. 

In Section~\ref{singularities} we deal with the elliptic obstructions
and we show, by an analysis of the toroidal compactifications, that
they disappear if $n\ge 9$, and also that there are no fixed divisors
at infinity. 

In Section~\ref{reflections} we examine the reflective obstructions by
describing the fixed divisors.  We do this first for arbitrary~$L$ and
then in greater detail for $L_{2d}$. 

In Section~\ref{specialcusp} we turn to the cusp obstructions. We
describe the structure of the cusps for a lattice $L$ having only
cyclic isotropic subgroups in its discriminant group. 

In Section~\ref{spinK3} we study the moduli space $\cS\cF_{2d}$ of
$\Kthree$ surfaces with a spin structure. In this case there are few
reflective obstructions, and the cusp forms constructed by Jacobi
lifting already have the properties we need. 

In Section~\ref{Borcherds} we show how to construct forms with the
properties needed for $\cF_{2d}$ by pulling back the Borcherds form. 
This requires us to find a suitable embedding of $L_{2d}$ in
$L_{2,26}$, which in turn requires a vector in $E_8$ with square $2d$
that is orthogonal to at most $12$ and at least $2$ roots.  We show
directly that such a vector exists for large~$d$ and use a small
amount of computer help to show that it exists for smaller~$d$. For
some values of $d$ we can find only a vector of square~$2d$ orthogonal
to $14$ roots. In these cases we can deduce that $\cF_{2d}$ has
non-negative Kodaira dimension.  \bigskip

\noindent{\it Acknowledgements:} We have learned much from
conversations with many people, but from S.~Kondo and
N.I.~Shepherd-Barron especially. We are grateful for financial support
from the Royal Society and the DFG Schwer\-punkt\-pro\-gramm SPP 1094
{\lq\lq{}Globale} Methoden in der komplexen {Geometrie{}\rq\rq}, Grant
Hu~337/5-3.  We are also grateful for the hospitality and good working
conditions provided by several places where one or more of us did
substantial work on this project: the Max-Planck-Institut f{\"{u}}r
Mathematik in Bonn; DPMMS in Cambridge and Trinity College, Cambridge;
Nagoya University; KIAS in Seoul; Tokyo University; and the Fields
Institute in Toronto. 

\section{Orthogonal groups and modular forms}\label{orthogonal}

Let $L$ be a lattice of signature $(2,n)$, with $n>1$. For any lattice
$M$ and field $K$ we write $M_K$ for $M\tensor K$.  Then $\cD_L$ is
one of the two connected components of
$$
\{[\bdw]\in \PP(L_\CC)\mid (\bdw,\bdw)=0,\ (\bdw,\Bar\bdw)>0\}. 
$$ 
We denote by $\Orth^+(L)$ the subgroup of $\Orth(L)$ that
preserves~$\cD_L$.  If $\Gamma<\Orth^+(L)$ is of finite index we
denote by $\cF_L(\Gamma)$ the quotient $\Gamma\backslash \cD_L$, which
is a quasi-projective variety by~\cite{BB}. 

For every non-degenerate integral lattice we denote by $L^\vee=\Hom
(L, \ZZ)$ its dual lattice. The finite group $A_L=L^\vee/L$ carries a
discriminant quadratic form $q_L$ (if $L$ is even) and a discriminant
bilinear form $b_L$, with values in $\QQ/2\ZZ$ and $\QQ/\ZZ$
respectively (see \cite [\S{}1.3]{Nik2}). We define
\begin{eqnarray*}
  \Tilde{\Orth}(L) &=& \{g\in \Orth(L)\mid g|_{A_L}=\id\},\ \text{ and}\\
  \Tilde\Orth^+(L) &=& \Tilde\Orth(L)\cap \Orth^+(L). 
\end{eqnarray*}
The $\Kthree$ lattice is
$$
L_{\Kthree}=3U \oplus 2E_8 (-1)
$$ 
where $U$ is the hyperbolic plane and $E_8$ is the (positive definite)
$E_8$-lattice. If $h\in L_{\Kthree}$ is a primitive vector with
$h^2=2d>0$ then its orthogonal complement $h^\perp_{L_{\Kthree}}$ is
isometric to
$$
L_{2d}=\latt{-2d} \oplus 2U \oplus 2E_8 (-1). 
$$
By \cite[Proposition 1.5.1]{Nik2}
$$
\Tilde{\Orth}(L_{2d})\imic\{ g\in \Orth(L_{\Kthree}) \mid g(h)=h \},
$$
and the moduli space $\cF_{2d}$ is given by
$$
\cF_{2d}=\Tilde\Orth^+(L_{2d}) \backslash \cD_{L_{2d}}. 
$$

A modular form of weight $k$ and character $\chi\colon \Gamma\to
\CC^*$ for a subgroup $\Gamma<\Orth^+(L)$ is a holomorphic function
$F\colon\cD_L^\bullet\to \CC$ on the affine cone $\cD_L^\bullet$ over
$\cD_L$ such that
\begin{equation}\label{mod-form}
  F(tZ)=t^{-k}F(Z)\ \forall\,t\in \CC^*,\ \text{ and }\
  F(gZ)=\chi(g)F(Z)\ \forall\,g\in \Gamma. 
\end{equation}
A modular form is a cusp form if it vanishes at every cusp. We denote
the linear spaces of modular and cusp forms of weight $k$ and
character $\chi$ for $\Gamma$ by $M_k(\Gamma,\chi)$ and
$S_k(\Gamma,\chi)$ respectively. 

\begin{theorem}\label{general_gt}
  Let $L$ be an integral lattice of signature $(2,n)$, $n\ge 9$, and
  let $\Gamma$ be a subgroup of finite index of $\Orth^+(L)$.  The
  modular variety $\cF_L(\Gamma)$ is of general type if there exists a
  character $\chi$ of finite order and a non-zero cusp form $F_a\in
  S_a(\Gamma,\chi)$ of weight $a<n$ that vanishes along the branch
  divisor of the projection $\pi\colon \cD_L\to \cF_L(\Gamma)$. 

  If $S_n(\Gamma,\det)\neq 0$ then the Kodaira dimension of
  $\cF_L(\Gamma)$ is non-negative. 
\end{theorem}
\begin{proof}
  We let $\Bar\cF_L(\Gamma)$ be a toroidal compactification of
  $\cF_L(\Gamma)$ with canonical singularities and no branch divisors
  at infinity, which exists by Theorem~\ref{main_sings_theorem}. We
  take a smooth projective model $\Hat\cF_L(\Gamma)$ by taking a
  resolution of singularities of $\Bar\cF_L(\Gamma)$. 

  Suppose that $F_{nk}\in M_{nk}(\Gamma, \det^k)$.  Then, if $dZ$ is a
  holomorphic volume element on $\cD_L$, the differential form
  $\Omega(F_{nk})=F_{nk}\,(dZ)^k$ is $\Gamma$-invariant and therefore
  determines a section of the pluricanonical bundle
  $kK=kK_{\Hat\cF_L(\Gamma)}$ away from the branch locus of
  $\pi\colon\cD_L\to \cF_L(\Gamma)$ and the cusps. 

  In general $\Omega(F_{nk})$ will not extend to a global section of
  $kK$.  We distinguish three kinds of obstruction to its doing so. 
  There are elliptic obstructions, arising because of singularities
  given by elliptic fixed points of the action of $\Gamma$; reflective
  obstructions, arising from the branch divisors in $\cD_L$ (divisors
  fixed pointwise by an element of $\Gamma$ acting locally as a
  quasi-reflection); and cusp obstructions, arising from divisors at
  infinity. 
 
  In this situation the elliptic obstruction vanishes (and there are
  no elliptic or reflective obstructions at infinity either) because
  of the choice of $\Bar\cF_L(\Gamma)$.  So $\Omega(F_{nk})$ will
  extend to a section of $kK$ provided it extends to a general point
  of each branch divisor and each boundary divisor. 

  We apply the low-weight cusp form trick, used for example in
  \cite{G}, \cite{GH1}, \cite{GS} to show that the cusp obstruction
  for continuation of the pluricanonical forms on a smooth
  compactification is small compared with the dimension of
  $S_{nk}(\Gamma, \det^k)$.  Let $N$ be the order of $\chi$ and put
  $k=2Nl$. Then we consider special elements $F^0_{nk}\in
  S_{nk}(\Gamma)$ of the form
\begin{equation}\label{verycusp}
  F^0_{nk}=F_a^k F_{(n-a)k}
\end{equation}
where $F_{(n-a)k}\in M_{(n-a)k}(\Gamma)$ is a modular form of weight
$(n-a)k\ge k$.  The corresponding differential form $\Omega(F^0_{nk})$
vanishes to order at least $k$ on the boundary of the toroidal
compactification $\overline\cF_L(\Gamma)$. It follows by the results
of \cite{AMRT} that $\Omega(F^0_{nk})$ extends as a $k$-fold
pluricanonical form to the generic point of any boundary divisor of
$\Bar\cF_L(\Gamma)$. The reason is that the anticanonical divisor of a
toric variety is the sum of the torus-invariant divisors, so $dZ$ has
simple poles at all boundary divisors in a toroidal compactification. 

Since $F_a$ vanishes at the branch divisors, which are the fixed
divisors of reflections by Theorem~\ref{qrefs_have_order_2},
$\Omega(F^0_{nk})$ vanishes there to order~$k$, and hence it extends
to give a section of $kK$ over $\Hat\cF_L(\Gamma)$. 

Finally, we observe that this gives us an injective map
$$
M_{(n-a)k}(\Gamma)\hookrightarrow H^0(\widehat\cF_L(\Gamma)). 
$$
But $\dim M_{(n-a)k}(\Gamma)\sim k^n$, as can be seen from \cite{BB}:
a more precise estimate, using the results of~\cite{Mum}, can be found
in~\cite{GHS1}. Hence it follows that $\cF_L(\Gamma)$ is of general
type.

Even if we can only find a cusp form of weight~$n$ we still get some
information, because of the well-known result of Freitag that if
$F_n\in S_n(\Gamma, \det)$ then $\Omega(F_n)$ defines an element of
$H^0(K_{\widehat\cF_L(\Gamma)})$. Therefore the plurigenera do not all
vanish: indeed $p_g\ge 1$. 
\end{proof}

\section{Singularities of locally symmetric varieties}\label{singularities}

In this section, we consider the singularities of compactified locally
symmetric varieties associated with the orthogonal group of a lattice
of signature~$(2,n)$. Our main theorem is that for all but small~$n$,
the compactification may be chosen to have canonical singularities. 

\begin{theorem}\label{main_sings_theorem}
  Let $L$ be a lattice of signature $(2,n)$ with $n\ge 9$, and let
  $\Gamma<\Orth^+(L)$ be a subgroup of finite index. Then there exists
  a toroidal compactification $\Bar\cF_L(\Gamma)$ of
  $\cF_L(\Gamma)=\Gamma\backslash\cD_L$ such that $\Bar\cF_L(\Gamma)$
  has canonical singularities and there are no branch divisors in the
  boundary. The branch divisors in $\cF_L(\Gamma)$ arise from the
  fixed divisors of reflections.
\end{theorem}
\begin{proof}
  Immediate from Corollaries~\ref{can_sings_on_interior},
  \ref{can_sings_on_dim0_cusps} and \ref{can_sings_on_dim1_cusps}. The
  last part is a summary of Theorem~\ref{qrefs_have_order_2} (an
  element that fixes a divisor in $\cD_L$ has order~$2$ on the tangent
  space) and Corollary~\ref{qrefs_come_from_refs_in_L} (such elements,
  up to sign, are given by reflections by vectors in~$L$). 
\end{proof}

In fact we prove more than this: for example, $\cF_L(\Gamma)$ has
canonical singularities if $n\ge 7$
(Corollary~\ref{can_sings_on_interior}), and our method (which uses
ideas from~\cite{Nik1}) gives some information about what
non-canonical singularities can occur for small~$n$. In order to
choose $\Bar\cF_L(\Gamma)$ as in Theorem~\ref{main_sings_theorem} it
is enough to take all the fans defining the toroidal compactification
to be basic.

\subsection{The interior}\label{interior}

For $[\bdw]\in \cD_L$ we define $\WW=\CC.\bdw$. We put
$S=(\WW\oplus\Bar\WW)^\perp\cap L$, noting that $S$ could be $\{0\}$,
and take $T=S^\perp\subset L$. 

In the case of polarised $\Kthree$ surfaces, $S$ is the primitive part
of the Picard lattice and $T$ is the transcendental lattice of the
surface corresponding to the period point~$\bdw$. 

\begin{lemma}\label{S_T_are_disjoint}
$S_\CC\cap T_\CC=\{0\}$. 
\end{lemma}

\begin{proof}
  $S_\CC$ and $T_\CC$ are real ({i.e.} preserved by complex
  conjugation) so it is enough to show that $S_\RR\cap T_\RR=\{0\}$. 
  If $\bdx\in T_\RR\cap S_\RR$ then $(\bdx,\bdx)=0$ from the
  definition of $T$, so it is enough to prove that $S_\RR$ is negative
  definite. The subspace $\UU=\WW\oplus\Bar\WW\subset L_\CC$ is also
  real, so we may write $\UU=U_\RR\tensor\CC$, taking $U_\RR$ to be
  the real vector subspace of $\UU$ fixed pointwise by complex
  conjugation. An $\RR$-basis for $U_\RR$ is given by
  $\{\bdw+\bar\bdw, i(\bdw-\bar\bdw)\}$.  But
  $(\bdw+\bar\bdw,\bdw+\bar\bdw)>0$ and $(i(\bdw-\bar\bdw),
  i(\bdw-\bar\bdw))>0$, so $U_\RR$ has signature $(2,0)$. Hence
  $U_\RR^\perp$ has signature $(0,n)$, but $S_\RR\subset U_\RR^\perp$
  so $S_\RR$ is negative definite. 
\end{proof}
We are interested first in the singularities that arise at fixed
points of the action of $\Gamma$ on $\cD_L$. Suppose then that
$\bdw\in L_\CC$ and let $G$ be the stabiliser of $[\bdw]$ in $\Gamma$. 
Then $G$ acts on $\WW$ and we let $G_0$ be the kernel of this action:
thus for $g\in G$ we have $g(\bdw)=\alpha(g)\bdw$ for some
homomorphism $\alpha\colon G \to \CC^*$, and $G_0=\ker \alpha$. 

\begin{lemma}\label{G_acts_on_T}
$G$ acts on $S$ and on $T$. 
\end{lemma}
\begin{proof}
  $G$ acts on $\WW$ and on $L$, hence also on
  $S=(\WW\oplus\Bar\WW)^\perp\cap L$ and on $T=S^\perp\cap L$. 
\end{proof}

\begin{lemma}\label{G0_acts_trivially_on_T}
  $G_0$ acts trivially on $T_\QQ$. 
\end{lemma}
\begin{proof}
  If $\bdx\in T_\QQ$ and $g\in G_0$ then
$$
(\bdw, \bdx)=(g(\bdw), g(\bdx))= (\bdw, g(\bdx)). 
$$
Hence $T_\QQ \ni \bdx-g(\bdx)\in L_\QQ\cap (\WW\oplus\Bar\WW) =S_\QQ$,
so by Lemma~\ref{S_T_are_disjoint} we have $g(\bdx)=\bdx$. 
\end{proof}

The quotient $G/G_0$ is a subgroup of $\Aut\WW\imic\CC^*$ and is thus
cyclic of some order, which we call $r_\bdw$. So by the above,
$\mu_{r_\bdw}\imic G/G_0$ acts on $T_\QQ$.  (By $\mu_r$ we mean the
group of $r$th roots of unity in $\CC$.) 

For any $r\in \NN$ there is a unique faithful irreducible
representation of $\mu_r$ over $\QQ$, which we call $\cV_r$. The
dimension of $\cV_r$ is $\varphi(r)$, where $\varphi$ is the Euler
$\varphi$ function and, by convention, $\varphi(1)=\varphi(2)=1$. The
eigenvalues of a generator of $\mu_r$ in this representation are
precisely the primitive $r$th roots of unity: $\cV_1$ is the
{$1$-dimensional} trivial representation. Note that $-\cV_d=\cV_d$ if
$d$ is even and $-\cV_d=\cV_{2d}$ if $d$ is odd. 

\begin{lemma}\label{splitting_of_T}
  As a $G/G_0$-module, $T_\QQ$ splits as a direct sum of irreducible
  representations $\cV_{r_\bdw}$. In particular, $\varphi(r_\bdw)|\dim
  T_\QQ$. 
\end{lemma}
\begin{proof}
  We must show that no nontrivial element of $G/G_0$ has $1$ as an
  eigenvalue on $T_\CC$. Suppose that $g\in G\setminus G_0$ (so
  $\alpha(g)\neq 1$) and that $g(\bdx) =\bdx$ for some $\bdx\in
  T_\CC$. Then
$$
( \bdw,\bdx) =( g(\bdw), g(\bdx)) = \alpha(g)( \bdw, \bdx),
$$
so $( \bdw,\bdx) =0$, so $\bdx\in S_\CC\cap T_\CC=0$. 
\end{proof}

\begin{corollary}\label{T_splits_for_g}
  If $g\in G$ and $\alpha(g)$ is of order $r$ (so $r|r_\bdw$), then
  $T_\QQ$ splits as a $g$-module into a direct sum of irreducible
  representations $\cV_r$ of dimension~$\varphi(r)$. 
\end{corollary}
\begin{proof}
  Identical to the proof of Lemma~\ref{splitting_of_T}. 
\end{proof}

We are interested in the action of $G$ on the tangent space to
$\cD_L$.  We have a natural isomorphism
$$
T_{[\bdw]}\cD_L\imic \Hom(\WW,\WW^\perp/\WW)=:V. 
$$
We choose $g\in G$ of order $m$ and put $\zeta=e^{2\pi i/m}$ for
convenience: as $g$ is arbitrary there is no loss of generality. Let
$r$ be the order of $\alpha(g)$, as in Corollary~\ref{T_splits_for_g}
(this is called $m$ in \cite{Nik1} but we want to keep the notation of
\cite{Ko1}). In particular $r|m$. The eigenvalues of $g$ on $V$ are
powers of $\zeta$, say $\zeta^{a_1},\ldots,\zeta^{a_n}$, with $0\le
a_i<m$. We define
\begin{equation}\label{RST}
  \Sigma(g):=\sum_{i=1}^n a_i/m. 
\end{equation}
Recall that an element of finite order in $\GL_n(\CC)$ (for any $n$)
is called a quasi-reflection if all but one of its eigenvalues are
equal to~$1$. It is called a reflection if the remaining eigenvalue is
equal to~$-1$. The branch divisors of $\cD_L\to\cF_L(\Gamma)$ are
precisely the fixed loci of elements of $\Gamma$ acting as
quasi-reflections. 
\begin{proposition}\label{RST_for_g_nonqref_and_phi(r)>4}
  Assume that $g\in G$ does not act as a quasi-reflection on $V$ and
  that $\varphi(r)>4$. Then $\Sigma(g)\ge 1$. 
\end{proposition}

\begin{proof}
  As $\xi$ runs through the $m$th roots of unity, $\xi^{m/r}$ runs
  through the $r$th roots of unity. We denote by
  $k_1,\ldots,k_{\varphi(r)}$ the integers such that $0<k_i<r$ and
  $(k_i,r)=1$, in no preferred order. Without loss of generality, we
  assume $\alpha(g)=\zeta^{mk_2/r}$ and
  $\Bar{\alpha(g)}=\alpha(g)^{-1}=\zeta^{mk_1/r}$, with $k_1\equiv
  -k_2 \bmod r$. 

  One of the $\QQ$-irreducible subrepresentations of $g$ on $L_\CC$
  contains the eigenvector $\bdw$: we call this $\VV_r^\bdw$ (it is
  the smallest $g$-invariant complex subspace of $L_\CC$ that is
  defined over $\QQ$ and contains $\bdw$). It is a copy of
  $\cV_r\tensor\CC$: to distinguish it from other irreducible
  subrepresentations of the same type we write
  $\VV_r^\bdw=\cV_r^\bdw\tensor\CC$. 

  If $\bdv$ is an eigenvector for $g$ with eigenvalue
  $\zeta^{mk_i/r}$, $i\neq 1$ (in particular $\bdv\not\in\Bar\WW$),
  then $\bdv\in \WW^\perp$ since
  $(\bdv,\bdw)=(g(\bdv),g(\bdw))=\zeta^{mk_i/r}\alpha(g)(\bdv, \bdw)$. 
  Therefore the eigenvalues of $g$ on $\VV_r^\bdw\cap \WW^\perp /\WW$
  include $\zeta^{mk_i/r}$ for $i\ge 3$, so the eigenvalues on
  $\Hom(\WW, \VV_r^\bdw\cap\WW^\perp/\WW)\subset V$ include
  $\zeta^{mk_1/r}\zeta^{mk_i/r}$ for $i\ge 3$. So, writing $\{a\}$ for
  the fractional part of $a$, we have
\begin{eqnarray}\label{W-contribution}
  \Sigma(g)&\ge&
  \sum_{i=3}^{\varphi(r)}\frac{1}{m}\left\{\frac{mk_1}{r} 
    +\frac{mk_i}{r}\right\}\nonumber\\
  &=&
  \sum_{i=3}^{\varphi(r)}\left\{\frac{k_1+k_i}{r}\right\}. 
\end{eqnarray}
Now the proposition follows from the elementary Lemma~\ref{bigphi}
below. 
\end{proof}

\begin{lemma}\label{bigphi}
  Suppose $k_1,\ldots,k_{\varphi(r)}$ are the integers between $0$ and
  $r$ coprime to $r$, in some order, and that $k_2=r-k_1$. If
  $\varphi(r)\ge 6$ then
$$
\sum_{i=3}^{\varphi(r)}\Big\{\frac{k_1}{r}+\frac{k_i}{r}\Big\}\ge 1. 
$$
\end{lemma}

\begin{proof}
  If $k_1<k_3<r/2$ then $\left\{\frac{k_1+k_3}{r}\right\}
  =\frac{k_1+k_3}{r}$, and $k_4=r-k_3$ so
  $\left\{\frac{k_1+k_4}{r}\right\} =\frac{k_1+r-k_3}{r}$. Thus
$$
\left\{\frac{k_1+k_3}{r}\right\}+\left\{\frac{k_1+k_4}{r}\right\}
=\frac{2k_1+r}{r}>1. 
$$
If $r/2>k_1>r/4$ or $r>k_1>3r/4$ then $(k_1+k_3)+(k_1+k_4)\equiv
2k_1\bmod r$, so
$$
\left\{\frac{k_1+k_3}{r}\right\}+\left\{\frac{k_1+k_4}{r}\right\}
\equiv\frac{2k_1}{r}\bmod 1. 
$$
Therefore
$\left\{\frac{k_1+k_3}{r}\right\}+\left\{\frac{k_1+k_4}{r}\right\}
>\frac12$, and similarly for
$\left\{\frac{k_1+k_5}{r}\right\}+\left\{\frac{k_1+k_6}{r}\right\}$,
so the sum is at least~$1$. 

If $r/2<k_1<3r/4$ then we may take $k_3=1$ and $k_4=r-1$, and then
$\left\{\frac{k_1+k_3}{r}\right\}+\left\{\frac{k_1+k_4}{r}\right\}
=1+\frac{2k_1}{r}>1$. 

The remaining possibility is that $k_1<r/4$ but $k_1>k_j$ if
$k_j<r/2$. But then there is no integer coprime to $r$ between $r/4$
and $3r/4$. As long as $2\lceil r/4\rceil<\lfloor 3r/4\rfloor$, which
is true if $r>9$, we may choose a prime $q$ such that $r/4<q<3r/4$, by
Bertrand's Postulate \cite[Theorem 418]{HW}, and $\gcd(q,r)\neq 1$ so
$r=2q$ or $r=3q$. In the first case one of $q\pm 2$ lies in
$(r/4,3r/4)$ and is prime to $r$, and in the second case one of $q\pm
1$ or $q\pm 2$ does, unless $r<8$; so this possibility does not occur. 
The cases $r=7$ and $r=9$, which are not covered by this argument, are
readily checked: $2\in (7/4,21/4)$ and $4\in(9/4,27/4)$ are coprime
to~$r$. 
\end{proof}

\begin{proposition}\label{RST_for_g_nonqref_and_r=1,2}
  Assume that $g\in G$ does not act as a quasi-reflection on $V$ and
  that $r=1$ or $r=2$. Then $\Sigma(g)\ge 1$. 
\end{proposition}

\begin{proof}
  We note first that we may assume $g$ is not of order~$2$, because if
  $g^2$ acts trivially on $V$ but $g$ is not a quasi-reflection then
  at least two of the eigenvalues of $g$ on $V$ are $-1$, and hence
  $\sum_{i=1}^n a_i/m\ge 1$. However, $g^2$ does act trivially
  on~$T_\CC$, by Corollary~\ref{T_splits_for_g}. Therefore $g^2$ does
  not act trivially on~$S_\CC$. The representation of $g$ on $S_\CC$
  therefore splits over $\QQ$ into a direct sum of irreducible
  subrepresentations $\cV_d$, and at least one such piece has $d>2$. 
  So on the subspace $\Hom(\WW,\cV_d\otimes
  \CC)=\Hom(\WW,(\cV_d\tensor\CC \oplus\WW)/\WW)\subset V$, the
  representation of $g$ is $\pm\cV_d$ (the sign depending on whether
  $r=1$ or $r=2$), and choosing two conjugate eigenvalues $\pm\zeta^a$
  and $\pm\zeta^{m-a}$ we have $\sum a_i/m\ge 1$. 
\end{proof}

\begin{theorem}\label{RST_for_nonqrefs_on_interior}
  Assume that $g\in G$ does not act as a quasi-reflection on $V$ and
  that $n\ge 6$. Then $\Sigma(g)\ge 1$. 
\end{theorem}
\begin{proof}
  In view of Proposition~\ref{RST_for_g_nonqref_and_r=1,2} and
  Proposition~\ref{RST_for_g_nonqref_and_phi(r)>4}, we need only
  consider $r=3$, $4$, $5$, $6$, $8$, $10$ or $12$. We suppose, as
  before, that $g$ has order $m$, and we put $k=m/r$. 

  Consider first a $\QQ$-irreducible subrepresentation $\cV_d\subset
  S_\CC$, and the action of $g$ on $\Hom(\WW,\cV_d\tensor \CC)\subset
  V$. This is $\zeta^{kc}\cV_d$, where $\zeta$ is a primitive $m$th
  root of unity, and $c$ is some integer with $0<c<r$ and $(c,r)=1$
  (the eigenvalue of $g$ on $\WW$ is $\zeta^{-kc}$.  So the
  eigenvalues are of the form $\zeta^{b_i/m}$ for $1\le i\le
  \varphi(d)$, with $0\le b_i<m$ and the $b_i$ all different mod~$m$
  but all equivalent mod~$l$, where $l=m/d$. Clearly
$$
\sum_{i=1}^{\varphi(d)}\frac{b_i}{m}\ge \frac{1}{2m}l(\varphi(d)-1)\varphi(d)
=\frac{1}{2d}(\varphi(d)-1) \varphi(d)
$$
and it is easy to see that this is $\ge 1$ unless
$d\in\{1,\ldots,6,8,10,12,18,30\}$. 

By a slightly less crude estimate we can reduce further. For $d>2$ we
write $c_{\min}(d)$ for a lower bound for the contribution to the sum
$\Sigma(g)$ from $\cV_d$ as a subrepresentation of $g$ on $S_\CC$,
{i.e.} 
$$
c_{\min}(d)=\min_{0\le a<d}\sum_{0<b<d,\ (d,b)=1}\left\{\frac{b+a}{d}\right\}. 
$$
Note that this is a lower bound independently of $r$. For fixed $r$
one has a contribution to $\Sigma(g)$ from $\cV_d$ of at most
\begin{eqnarray*}
  \min_{0<c<r}\sum_{0<b<d,\ (d,b)=1}\left\{\frac{bl+kc}{m}\right\}
  &=&\min_{0<c<r}\sum_{0<b<d,\ (d,b)=1}\left\{\frac{b}{d}
    +\frac{kc}{m}\right\}\\
  &\ge& \min_{0<c<r}\sum_{0<b<d,\
    (d,b)=1}\left\{\frac{b}{d}+\frac{\lfloor kc/l \rfloor}{d}\right\}\\
  &\ge& c_{\min}(d). 
\end{eqnarray*}

It is easy to calculate that $c_{\min}(30)=92/30$ (attained when
$a=19$), $c_{\min}(18)=42/18$, $c_{\min}(12)=16/12$,
$c_{\min}(10)=12/10$, $c_{\min}(8)=12/8$ and $c_{\min}(5)=6/5$. But
\begin{equation}\label{cmins}
  c_{\min}(3)=c_{\min}(6)=1/3, \quad c_{\min}(4)=1/2. 
\end{equation}
Hence we may assume that $r\in\{3,4,5,6,8,12\}$ and
$d\in\{1,2,3,4,6\}$ for every subrepresentation $\cV\tensor\CC\subset
S_\CC$. The summands of $T_\CC$ are all $\cV_r\tensor\CC$. We let
$\nu_d$ be the multiplicity of $\cV_d$ in $S_\CC$ as a $g$-module, and
$\lambda$ be the multiplicity of $\cV_r$ in $T_\CC$. Counting
dimensions gives
\begin{equation}\label{dimcount}
  \lambda\varphi(r)+\nu_1+\nu_2+2\nu_3+2\nu_4+2\nu_6=n+2. 
\end{equation}
We split into two cases, depending on whether $\varphi(r)=4$ or
$\varphi(r)=2$.  
\medskip

\noindent{\bf Case I.} Suppose $\varphi(r)=4$, so $r\in\{5,8,10,12\}$. 
\smallskip

If $\lambda>1$ then there will be a $\cV_r\tensor\CC$ not containing
$\WW$ and this will contribute at least $c_{\min}(r)$ to $\Sigma(g)$,
just as if it were contained in $S_\CC$ instead of $T_\CC$. For $r=5$,
$8$, $10$ or $12$ we have $c_{\min}(r)\ge 1$, so we may assume that
$\lambda=1$.  Moreover in these cases $\varphi(r)=4$, so
equation~(\ref{dimcount}) becomes
\begin{equation}\label{dimcount_phi=4}
  \nu_1+\nu_2+2\nu_3+2\nu_4+2\nu_6=n-2. 
\end{equation}
We may assume that $\nu_4\le 1$ and $\nu_3+\nu_6\le 2$, as otherwise
those summands contribute at least~$1$ to $\Sigma(g)$, by
equation~(\ref{cmins}). The contribution from $\cV^\bdw_r$ was
computed in equation~(\ref{W-contribution}) above: for $\varphi(r)=4$
it is $\frac{k_1+k_3}{r}+\frac{k_1+k_4}{r}$. The contribution from a
$\cV_1$ (an invariant) is $\frac{k_1}{r}$ and from $\cV_2$ (an
anti-invariant) it is $\{\frac{k_1}{r}+\frac{1}{2}\}$. 

Now we can compute all cases. The contribution from a copy of $\cV_d$
is
\begin{equation}\label{d-contribution}
  \sum_{(a,d)=1}\left\{\frac{a}{d}+\frac{k_1}{r}\right\}
\end{equation}
or $\frac{k_1}{r}$ if $d=1$. Half the time ($k_1$ first or third in
order of size) the contribution $c^\bdw$ from $\cV^\bdw_r$ is already
at least~$1$. In all cases it is at least $\frac{1}{2}$, so we may
also assume that $\nu_4=0$. In six of the remaining eight cases we get
$\Sigma(g)\ge 1$ unless $L_\CC=\VV_r^\bdw$ and hence $n=2$: all other
possible contributions are greater than $1-c^\bdw$. The exceptions are
$r=5$, $k_1=4$ and $r=10$, $k_1=3$. 

For $r=5$, $k_1=4$, contributions from $\cV_r^\bdw$, $\cV_1$, $\cV_2$,
$\cV_3$ and $\cV_6$ are $\frac{3}{5}$, $\frac{4}{5}$, $\frac{3}{10}$,
$\frac{3}{5}$ and $\frac{8}{5}$ respectively. So $\Sigma(g)\ge 1$
unless $\nu_1=\nu_3=\nu_6=0$ and $\nu_2\le 1$, and in particular $n\le
3$. 

For $r=10$, $k_1=3$, contributions from $\cV_r^\bdw$, $\cV_1$,
$\cV_2$, $\cV_3$ and $\cV_6$ are $\frac{3}{5}$, $\frac{3}{10}$,
$\frac{8}{10}$, $\frac{6}{10}$ and $\frac{6}{10}$ respectively. So
$\Sigma(g)\ge 1$ unless $\nu_2=\nu_3=\nu_6=0$ and $\nu_1\le 1$, and in
particular $n\le 3$.  \medskip

\noindent{\bf Case II.} Suppose $\varphi(r)=2$, so $r\in\{3,4,6\}$. 
\smallskip

In this case one summand of $L_\CC$ as a $g$-module is the space
$\WW\oplus\Bar\WW$, which is $\VV_r^\bdw$, a copy of
$\cV_r\tensor\CC$. We denote by $\nu_d$ the multiplicity of $\cV_d$ in
$L_\CC/(\WW\oplus\Bar\WW)$ as a $g$-module. Thus $\nu_r$ is the number
of copies of $\cV_r\tensor\CC$ in $L_\CC$ that are different from
$\VV_r^\bdw$. Equation~(\ref{dimcount}) becomes
\begin{equation}\label{dimcount_phi=2}
\nu_1+\nu_2+2\nu_3+2\nu_4+2\nu_6=n. 
\end{equation}
There are six cases (three values of $r$, and $k_1=1$ or $k_1=r-1$)
and we simply compute all contributions in each case using the
expression~(\ref{d-contribution}). For $1$-dimensional summands ($d=1$
or~$2$) the lowest contribution is $\frac{1}{6}$ (for $r=3$, $k_1=2$,
$d=2$ and for $r=6$, $k_1=1$ and $d=1$). For $2$-dimensional summands
the lowest contribution is $\frac{1}{3}$ (for $r=3$, $k_1=2$, $d=3$
and for $r=6$, $k_1=1$, $d=6$). So $\Sigma(g)\ge 1$ unless $n\le 5$. 
\end{proof}

\begin{corollary}\label{can_sings_off_branch_divisors}
  If $n\ge 6$, then the space $\cF_L(\Gamma)$ has canonical
  singularities away from the branch divisors of $\cD_L\to
  \cF_L(\Gamma)$. 
\end{corollary}
\begin{proof}
  This follows at once from the Reid--Shepherd-Barron--Tai criterion
  (RST criterion for short) for canonical singularities: see~\cite{Re}
  or~\cite{T}. 
\end{proof}

\noindent {\it Remark.} It is easy to classify the types of canonical
singularities that can occur for small~$n$, by examining the
calculations above. 

So far we have not considered quasi-reflections. We need to analyse
not only quasi-relections themselves but also all elements some power
of which acts as a quasi-reflection on~$V$: note, however, that
Theorem~\ref{RST_for_nonqrefs_on_interior} does apply to such
elements. 

\begin{theorem}\label{qrefs_have_order_2}
  Suppose $n>2$. Let $g\in G$ and suppose that $h=g^k$ acts as a
  quasi-reflection on~$V$. Then, as a $g$-module, $L_\QQ$ is either
  $\cV_k\oplus\bigoplus_j \cV_{2k}$ or $\cV_{2k}\oplus\bigoplus_j
  \cV_{k}$ (that is, one copy of $\cV_k$ and some copies of $\cV_{2k}$
  or vice versa). In particular, $h$ has order~$2$. 
\end{theorem}
\begin{proof}
  Suppose that $L_\QQ$ decomposes as a $g$-module as
  $\cV_r^\bdw\oplus\bigoplus_i \cV_{d_i}$ for some sequence
  $d_i\in\NN$.  The eigenvalues of $h$ on $V$ are all equal to~$1$,
  with exactly one exception. On the other hand, if $\zeta_r$ and
  $\zeta_{d_i}$ denote primitive $r$th and $d_i$th roots of unity, the
  eigenvalues of $h$ are certain powers of $\zeta_r$ (on
  $\Hom(\WW,\VV_r^\bdw\cap\WW^\perp/\WW)$) and all numbers of the form
  $\alpha(h)^{-1}\zeta_{d_i}^{ka}$ for $(a,d_i)=1$. 

  Consider a $\cV_d=\cV_{d_i}$ and put $d'=d/(k,d)$. The eigenvalues
  of $h$ on $\cV_d$ are primitive $d'$th roots of unity: each one
  occurs with multiplicity exactly $\varphi(d)/\varphi(d')$. However,
  only two eigenvalues of $h$ may occur in any $\cV_d$, and only one
  (namely $\alpha(h)$) may occur with multiplicity greater than~$1$,
  since if $\xi$ is an eigenvalue of~$h$ on $\cV_d$, the eigenvalue
  $\alpha(h)^{-1}\xi$ occurs with the same multiplicity on~$V$. Hence
  $\varphi(d')\le 2$, and if $\varphi(d')=2$ then $\varphi(d)=2$: this
  last can occur at most once. 

  Let us consider first the case where for some $d$ we have
  $\varphi(d)=\varphi(d')=2$. We claim that in this case $n=2$. We
  must have $d=6$ and $(k,d)=2$, and therefore $\alpha(h)=\omega$, a
  primitive cube root of unity. There can be no other $\cV_d$ summands
  ({i.e.} summands not containing $\WW$), because such a $\cV_d$ would
  have $\varphi(d)=1$ and hence give rise to an eigenvalue
  $\pm\omega^2$ for $h$ on $V$; but the $\cV_6$ already gives rise to
  an eigenvalue for $h$ on $V$ different from~$1$. So
  $L_\QQ=\cV_r^\bdw\oplus\cV_6$. The eigenvalues of $h$ on
  $\cV_r^\bdw$ are $\omega$ and $\omega^2$, each with multiplicity
  $\varphi(r)/2$: so $\varphi(r)=2$, otherwise $h$ has the eigenvalue
  $\omega$ with multiplicity~$>1$ on~$V$. Hence $\rank L=4$ and $n=2$. 

  Since we are assuming that $n\ge 6$, we have $\varphi(d')=1$ for
  all~$d$: that is, the eigenvalues of $h$ on the $\cV_d$ part are
  all~$\pm 1$. Put $r'=r/(k,r)$. We claim that $\varphi(r')=1$. 

  Suppose instead that $\varphi(r')\ge 2$, so $\alpha(h)\neq \pm 1$. 
  Then $\varphi(r)/\varphi(r')\le 2$, since the multiplicity of
  $\alpha(h)^{-2}\neq 1$ as an eigenvalue of $h$ on~$V$ is at least
  $\varphi(r)/\varphi(r')-1$. But the eigenvalues of $h$ on
  $\cV_r^\bdw$ are the primitive $r'$th roots of unity. If
  $\varphi(r')>2$ then these include $\alpha(h)$, $\alpha(h)^{-1}$,
  $\xi$ and $\xi^{-1}$ for some $\xi$, these being distinct. But then
  the eigenvalues of $h$ on $V$ include $\alpha(h)^{-1}\xi$ and
  $\alpha(h)^{-1}\xi^{-1}$, neither of which is equal to~$1$. So
  $\varphi(r')\le 2$

  Moreover, if $\varphi(r)/\varphi(r')=2$ then $h$ has the eigenvalue
  $\alpha(h)^{-2}\neq 1$ on~$V$, and any $\cV_d$ will give rise to the
  eigenvalue $\pm\alpha(h)^{-1}\neq 1$; so no such components occur,
  and $L_\QQ=\cV_r^\bdw$. Moreover, $\varphi(r)\le 4$ so $n\le 2$. 

  This shows that if $h$ is a quasi-reflection and $\varphi(r')>1$
  then $\varphi(r')=2$; moreover if $n>2$ then
  $\varphi(r)=\varphi(r')=2$.  Hence, if $\varphi(r')>1$, we have
  $r=6$ and $(r,k)=2$, so again $\alpha(h)=\omega$, a primitive cube
  root of unity. This time $\WW\oplus\Bar\WW=\VV_r^\bdw$, so the
  eigenvalues of $h$ on $V$ all arise from $\cV_d$ and since
  $\varphi(d')=1$ they are equal to $\pm\omega^2\neq 1$. So there is
  only one of them, that is, $n=1$. 

  Since we suppose $n>2$, it follows that $\varphi(r')=1$. The theorem
  follows immediately from this. 
\end{proof}

\begin{corollary}\label{qrefs_come_from_refs_in_L}
  The quasi-reflections on $V$, and hence the branch divisors of
  $\cD_L\to\cD_F(\Gamma)$, are induced by elements $h\in \Orth(L)$
  such that $\pm h$ is a reflection with respect to a vector in $L$. 
\end{corollary}

\begin{proof}
  The two cases are distinguished by whether $\alpha(h)=\pm 1$. If
  $\alpha(h)=1$ then the eigenvalues of $h$ on $L_\CC$ are $+1$ with
  multiplicity~$1$ and $-1$ with multiplicity~$n+1$, so $-h$ is a
  reflection; if $\alpha(h)=-1$, they are the other way round. 
\end{proof}

Now suppose that $g\in G$ and that $g^k=h$ is a quasi-reflection,
$k>1$.  By Theorem~\ref{qrefs_have_order_2}, $h$ has order~$2$ so $g$
has order $2k$. We may suppose that the eigenvalues of $g$ on~$V$ are
$\zeta^{a_1},\ldots,\zeta^{a_n}$, where $\zeta$ is a primitive $2k$th
root of unity, $0\le a_i<2k$, $a_n$ is odd and $a_i$ is even for
$i<n$. 

We need to look at the action of the group $\latt{g} /\latt{h}$ on
$V':=V/\latt{h}$. The eigenvalues of $g^l\latt{h}$ on $V'$ are
$\zeta^{la_1},\ldots,\zeta^{la_{n-1}}, \zeta^{2la_n}$, and we define
\begin{equation}\label{Sigmadash}
  \Sigma'(g^l):=\left\{\frac{la_n}{k}\right\}+
  \sum_{i=1}^{n-1}\left\{\frac{la_i}{2k}\right\}. 
\end{equation}
\begin{lemma}\label{modified_RST}
  $\cF_L(\Gamma)$ has canonical singularities if $\Sigma(g)\ge 1$ for
  every $g\in\Gamma$ no power of which is a quasi-reflection, and
  $\Sigma'(g^l)\ge 1$ if $g^k=h$ is a quasi-reflection and $1\le l<k$. 
\end{lemma}
\begin{proof}
  It is easy to see that if $V/\latt{g}$ has canonical singularities
  for every $g\in G$ then $V/G$ has canonical singularities (the
  converse is false). This follows from the fact that a $G$-invariant
  form extends to a resolution of $V/G$ if and only if it extends to a
  resolution of every $V/\latt{g}$, which is
  \cite[Proposition~3.1]{T}. 

  If no power of $g$ is a quasi-reflection on $V$ we simply apply the
  RST criterion. Otherwise, consider $g$ with $g^k=h$ a
  quasi-reflection as above. By
  Corollary~\ref{qrefs_come_from_refs_in_L}, $V'$ is smooth, and
  $V/\latt{g}\imic V'/(\latt{g}/\latt{h})$. So the result follows by
  applying the RST criterion to the elements $g^l\latt{h}$ acting
  on~$V'$. 
\end{proof}

\begin{proposition}\label{g^l_satisfies_modified_RST}
  If $g^k=h$ is a quasi-reflection and $n\ge 7$ then $\Sigma'(g^l)\ge
  1$ for every $1\le l<k$. 
\end{proposition}
\begin{proof}
  In fact we shall show that $\sum_{i=1}^{n-1}\{\frac{la_i}{2k}\}\ge
  1$.  As in Corollary~\ref{qrefs_come_from_refs_in_L} we have
  $\alpha(h)=\pm 1$ and this is a primitive $r'$th root of unity; so
  all the eigenvalues of $h$ on $\cV_r^\bdw$ are equal to $\alpha(h)$. 
  Here, as usual, $\WW\oplus\Bar\WW\subset\VV_r^\bdw$ (two copies of
  $\cV_r\otimes\CC$ if $r|2$) and we have decomposed $L_\CC$ as a
  $g$-module into $\QQ$-irreducible pieces. But exactly one eigenvalue
  of $h$ on $L_\CC$ is $-\alpha(h)=\mp 1$, and this must occur on some
  summand $\cV_d$. 

  The eigenvalues of $g$ on $\cV_d$ are primitive $d$th roots of
  unity, and in particular they all have the same order. Therefore the
  eigenvalues of $h$ are either all equal to~$1$ (if $\alpha(h)=-1$
  and $d|k$) or all equal to~$-1$ (if $\alpha(h)=1$ and $d|2k$ but $d$
  does not divide $k$). Since the eigenvalue $-\alpha(h)$ on $L_\CC$
  has multiplicity~$1$, it follows that $\varphi(d)=1$, {i.e.} $d=1$
  or $d=2$. 

  The eigenvector in $V$ corresponding to $\zeta^{a_n}$ comes from
  $\cV_d$, {i.e.} its span is the space
  $\Hom(\WW,\cV_d\tensor\CC)\subset V$. If we choose a primitive
  generator $\delta$ of $\cV_d\cap L$ we have $\delta^2<0$ since
  $\cV_d\subset U_\QQ^\perp$ as in Lemma~\ref{S_T_are_disjoint}, so
  $L'=\delta^\perp$ is of signature $(2,n-1)$ and $\latt{g}/\latt{h}$
  acts on $L'$ as a subgroup of $\Orth^+(L')$. But then
  $\Sigma'(g^l)=\{\frac{la_n}{k}\}+\Sigma(g^l\latt{h})$ where
  $g^l\latt{h}\in\Orth^+(L')$. It is clear that $g^l\latt{h}$ cannot
  be a quasi-reflection on $L'$: if it were, then by
  Corollary~\ref{qrefs_come_from_refs_in_L} the eigenvalues of $g^l$
  on $L'$ are all $\pm 1$, and so is its eigenvalue on $\cV_d$, so it
  has order dividing~$2$; so $g^l\in \latt{h}$. 

  Now we apply Theorem~\ref{RST_for_nonqrefs_on_interior} to $L'$,
  using $n-1\ge 6$. 
\end{proof}
\begin{corollary}\label{can_sings_on_interior}
  If $n\ge 7$ then $\cF_L(\Gamma)$ has canonical singularities. 
\end{corollary}

\subsection{Dimension~$0$ cusps}

We now consider the boundary
$\Bar\cF_L(\Gamma)\setminus\cF_L(\Gamma)$. Boundary components in the
Baily-Borel compactification correspond to totally isotropic subspaces
$E\subset L_\QQ$. Since $L$ has signature $(2,n)$, the dimension of
$E$ is $1$ or $2$, corresponding to dimension~$0$ and dimension~$1$
boundary components respectively. In this section we consider the case
$\dim E=1$, that is, isotropic vectors in~$L$. 

For a cusp~$F$ (of any dimension) we denote by $U(F)$ the unipotent
radical of the stabiliser subgroup $N(F)\subset \Gamma_\RR$ and by
$W(F)$ its centre. We let $N(F)_\CC$ and $U(F)_\CC$ be the
complexifications and put $N(F)_\ZZ=N(F)\cap \Gamma$ and
$U(F)_\ZZ=U(F)\cap \Gamma$. 

A toroidal compactification over a $0$-dimensional cusp $F$ coming
from a $1$-dimensional isotropic subspace~$E$ corresponds to an
admissible fan $\Sigma$ in some cone $C(F)\subset U(F)$. We have, as
in~\cite{AMRT}
$$
\cD_L(F)=U(F)_\CC\cD_L\subset \check\cD_L
$$
and in this case
$$
\cD_L(F)\cong F\times U(F)_\CC=U(F)_\CC. 
$$ 
Put $M(F)=U(F)_\ZZ$ and define the torus $\bdT(F)=U(F)_\CC/M(F)$. In
general $(\cD_L/M(F))_\Sigma$ is by definition the interior of the
closure of $\cD_L/M(F)$ in $\cD_L(F)/M(F)\cross_{\bdT(F)}
X_\Sigma(F)$, {i.e.} in $X_\Sigma(F)$ in this case, where
$X_\Sigma(F)$ is the torus embedding corresponding to the
torus~$\bdT(F)$ and the fan~$\Sigma$. We may choose $\Sigma$ so that
$X_\Sigma(F)$ is smooth and $G(F):=N(F)_\ZZ/U(F)_\ZZ$ acts on
$(\cD_L/M(F))_\Sigma$. The toroidal compactification is locally
isomorphic to $X_\Sigma(F)/G(F)$. Thus the problem of determining the
singularities is reduced to a question about toric varieties. The
result we want will follow from
Theorem~\ref{can_sings_on_toric_quotient}, below. We also need to
consider possible fixed divisors in the boundary. 

We take a lattice $M$ of dimension $n$ and denote its dual lattice by
$N$. A fan $\Sigma$ in $N\tensor\RR$ determines a toric variety
$X_\Sigma$ with torus $\bdT=\Hom(M,\CC^*)=N\tensor\CC^*$. 

\begin{theorem}\label{can_sings_on_toric_quotient}
  Let $X_\Sigma$ be a smooth toric variety and suppose that a finite
  group $G<\Aut(\bdT)=\GL(M)$ of torus automorphisms acts
  on~$X_\Sigma$.  Then $X_\Sigma/G$ has canonical singularities. 
\end{theorem}
\begin{proof}
  It is enough to show that for each $x\in X_\Sigma$ and for each
  $g\in \Stab_G(x)$, the quotient $X_\Sigma/\latt{g}$ has canonical
  singularities at~$x$. 

  We consider the subtorus $\bdT_0=\Stab_\bdT(x)$ of $\bdT$, which is
  given by $\bdT_0=N_0\tensor \CC^*$ for some sublattice $N_0\subset
  N$, and the quotient torus $\bdT_1=\bdT/\bdT_0$. The orbit
  $\orb(x)=\bdT.x$ of $x$ is isomorphic to $\bdT_1$: it corresponds to
  a cone $\sigma\in\Sigma$ of dimension
$$
s=\dim\sigma=\dim \bdT_1 = \codim\orb(x),
$$ 
and $N_0$ is the lattice generated by $\sigma\cap N$.  More
explicitly, $\orb(x)$ is given locally near $x$ by the equations
$\xi_i=0$, where $\xi_i$ are coordinates on~$\bdT_0$. The quotient
torus $\bdT_1$ is naturally isomorphic to $N_1\tensor \CC^*$, where
$N_1=N/N_0$ which is a lattice because $X_\Sigma$ is smooth. 

Certainly $x$ determines $\orb(x)$ and therefore $\sigma$, so $g$
stabilises~$\sigma$. If $U_\sigma=\Hom(M\cap \check\sigma,\CC^*)$
(semigroup homomorphisms) is the corresponding $\bdT$-invariant open
set, then $U_\sigma$ is $g$-invariant and the tangent spaces to
$U_\sigma$ and to $X_\Sigma$ at $x$ are the same: we denote this
tangent space by~$V$. Choosing a basis for $N_0$ and extending it to a
basis for $N$ gives an isomorphism of $U_\sigma$ with $\CC^s\times
(\CC^*)^{n-s}$ (compare~\cite[Theorem~1.1.10]{Od}). Since $g$
preserves $N_0$ it acts on both factors, by permuting the coordinates
and by torus automorphisms respectively. Thus
$$
V=(N_0\tensor\CC)\oplus\Lie(\bdT_1)
=(N_0\tensor\CC)\oplus(N_1\tensor\CC)=V_0\oplus V_1
$$
as a $g$-module, which is thus defined over~$\QQ$. 

Since $V$ is defined over $\QQ$, we may decompose it as a direct sum
of $\cV_d$s as a $g$-module, with each $d$ dividing $m$, the order
of~$g$. 

Note that if $g$ acts as a quasi-reflection, with eigenvalues
$(1,\ldots,1,\zeta)$ then since $g\in\GL(N)=\GL_n(\ZZ)$ we have
$\trace(g)=\zeta+n-1\in\ZZ$, and therefore $\zeta=-1$ and $g$ is a
reflection. 

We define $\Sigma(g)$ as we did in equation~(\ref{RST}) above, and in
the event that some power of $g$, say $h=g^k$, acts as a
quasi-reflection we define $V'=V/\latt{h}$ and $\Sigma'(g^l)$ as we
did in equation~(\ref{Sigmadash}). Now the theorem follows from
Proposition~\ref{toric_RST} and Proposition~\ref{toric_modified_RST},
below. 
\end{proof}

Note that we only needed to choose $\Sigma$ smooth: no further
subdivision is necessary. 

A version of Theorem~\ref{can_sings_on_toric_quotient} is stated
in~\cite{S-B} and proved in~\cite{Sn}. There the variety $X_\Sigma$ is
itself allowed to have canonical singularities, but $G$ is assumed to
act freely in codimension~$1$. 

\begin{proposition}\label{toric_RST}
  If $g\in G$ is not the identity, then unless $g$ acts as a
  reflection, $\Sigma(g)\ge 1$. 
\end{proposition}
\begin{proof}
  If $V$ contains a $\cV_d$ with $\varphi(d)>1$ then $g$ has a
  conjugate pair of eigenvalues and they contribute $1$ to
  $\Sigma(g)$. The same is true if $V$ contains two copies of $\cV_2$. 
  If neither of these is true, then $V=\cV_2\oplus(n-1)\cV_1$ and $g$
  is a reflection. 
\end{proof}

\begin{lemma}\label{toric_modified_RST}
  If $g^k=h$ acts as a reflection, and $g$ has order $m=2k>2$, then
  $\Sigma'(g^l)\ge 1$ for $1\le l<k$. 
\end{lemma}

\begin{proof}
  Since $m>2$, certainly $V$ contains a $\cV_d$ with $\varphi(d)\ge
  2$. In such a summand, the eigenvalues of any power of $g$ come in
  conjugate pairs: in particular, this is true for the eigenvalues of
  $h$. Therefore the eigenvalues of $h$ on $\cV_d$ are equal to~$1$ if
  $\varphi(d)\ge 2$, since the eigenvalue $-1$ occurs with
  multiplicity~$1$. Therefore a pair of conjugate eigenvalues of $g^l$
  on $\cV_d$ contribute~$1$ to $\Sigma'(g^l)$. 
\end{proof}

\begin{lemma}\label{no_fixed_toric_divisors}
  Let $X_\Sigma$ and $g$ be as above. Then there is no divisor in the
  boundary $X\setminus\bdT$ that is fixed pointwise by a non-trivial
  element of~$\latt{g}$. 
\end{lemma}
\begin{proof}
  Suppose $D$ were such a divisor, fixed pointwise by some element
  $h\in G$. Then $D$ corresponds to a $1$-parameter subgroup
  $\lambda\colon\CC^*\to \bdT$. Moreover, $D$ is a toric divisor and
  is itself a toric variety with dense torus $\bdT/\lambda(\CC^*)$. 
  Thus $h\in\GL(M)\imic \GL_n(\ZZ)$ acts trivially on
  $\bdT/\lambda(\CC^*)$; but the only such element is
  $\lambda(t)\mapsto \lambda(t^{-1})$, which does not preserve~$D$. 
\end{proof}

\begin{corollary}\label{can_sings_on_dim0_cusps}
  The toroidal compactification $\Bar\cF_L(\Gamma)$ may be chosen so
  that on a boundary component over a dimension~$0$ cusp,
  $\Bar\cF_L(\Gamma)$ has canonical singularities, and there are no
  fixed divisors in the boundary. 
\end{corollary}
\begin{proof}
  Since $\Sigma$ is $G(F)$-invariant, the result follows immediately
  from Theorem~\ref{can_sings_on_toric_quotient} and
  Lemma~\ref{no_fixed_toric_divisors}
\end{proof}

\begin{corollary}\label{no_fixed_divisors_at_dim0}
  There are no divisors at the boundary over a dimension~$0$ cusp~$F$
  that are fixed by a nontrivial element of~$G(F)$. 
\end{corollary}

Note that in this subsection we needed no restriction on~$n$. 

\subsection{Dimension~$1$ cusps}

It remains to consider the dimension~$1$ cusps. Here we have to be
more explicit: we consider a rank~$2$ totally isotropic subspace
$E_\QQ\subset L_\QQ$, corresponding to a dimension~$1$ boundary
component $F$ of $\cD_L$. We want to choose standard bases for $L_\QQ$
so as to be able to identify $U(F)$, $U(F)_\ZZ$ and $N(F)_\ZZ$
explicitly, as is done in~\cite{Sc} for maximal $\Kthree$ lattices,
where $n=19$. But we shall not be able to choose suitable bases of $L$
itself, as in~\cite{Sc}. The first steps, however, can be done
over~$\ZZ$. We define $E=E_\QQ\cap L$ and $E^\perp=E^\perp_\QQ\cap L$,
primitive sublattices of~$L$. 

\begin{lemma}\label{special_basis_for_L}
  There exists a basis $\bde_1',\ldots,\bde_{n+2}'$ for $L$ over $\ZZ$
  such that $\bde_1',\bde_2'$ is a basis for $E$ and
  $\bde_1',\ldots,\bde_n'$ is a basis for $E^\perp$. Furthermore we
  can choose $\bde_1',\ldots,\bde_{n+2}'$ so that
$$
A=\begin{pmatrix}\delta&0\\ 0&\delta e\end{pmatrix}
$$
for some integers $\delta$ and $e$, where $A$ is defined by
$$
Q':=(\bde_i',\bde_j')=\begin{pmatrix}0&0&A\\ 0&B&C\\
{}^t A&{}^tC&D\end{pmatrix}. 
$$
\end{lemma}
\begin{proof}
  We can find a basis with all the properties except for the special
  form of~$A$ by choosing any bases for the primitive sublattices $E$
  and $E^\perp$ of~$L$. Then the matrix $A$ may be chosen to have the
  special form given by choosing $\bde_1'$, $\bde_2'$, $\bde_{n+1}'$
  and $\bde_{n+2}'$ suitably: the numbers $\delta$ and $\delta e$ are
  the elementary divisors of $A\in\Mat_{2\times 2}(\ZZ)$. 
\end{proof}
If we are willing to allow two of the basis vectors to be in $L_\QQ$
we can achieve much more. 
\begin{lemma}\label{good_Q_basis}
  There is a basis $\bde_1,\ldots,\bde_{n+2}$ for $L_\QQ$ such that
  $\bde_1$ and $\bde_2$ form a $\ZZ$-basis for $E$, and
  $\bde_1,\ldots,\bde_n$ form a $\ZZ$-basis for $E^\perp$, for which
$$
Q:=(\bde_i,\bde_j)=\begin{pmatrix}
0&0&A\\ 0&B&0\\
A&0&0		 \end{pmatrix}
$$
with $A$ and $B$ as before. 
\end{lemma}
\begin{proof}
  We start with the basis $\bde_1',\ldots,\bde_{n+2}'$ from
  Lemma~\ref{special_basis_for_L}. Note that $B\in\Mat_{n-2 \times
    n-2}$ has non-zero determinant, because it represents the
  quadratic form of~$L$ on $E_\QQ^\perp/E_\QQ$. So we put
  $R=-B^{-1}C\in\Mat_{n-2\times 2}(\QQ)$ and we take $\bde_i$
  consisting of the columns of
$$
N:=\begin{pmatrix}I&0&R'\\ 0&I&R\\ 0&0&I\end{pmatrix},
$$
where $R'$ is chosen to satisfy
$$
D-{}^tCB^{-1}C+{}^tR'A+{}^tAR'=0. 
$$
Then $\bde_i$ is a $\QQ$-basis for $L_\QQ$ including $\ZZ$-bases for
$E$ and $E^\perp$, as we want, and ${}^tNQ'N=Q$ as required. 
\end{proof}

\begin{lemma}\label{describe_NUW_for_dim_1}
  The subgroups $N(F)$, $W(F)$ and $U(F)$ are given by
$$
N(F)=\left\{\begin{pmatrix}U&V&W\\ 0&X&Y\\ 0&0&Z\end{pmatrix}\mid
\begin{matrix}{}^t UAZ=A, {}^tXBX=B, {}^tXBY+{}^tVAZ=0,\\
{}^tYBY+{}^tZAW+{}^tWAZ=0,\ \det U>0\end{matrix}\right\},
$$
$$
W(F)=\left\{\begin{pmatrix}I&V&W\\ 0&I&Y\\ 0&0&I\end{pmatrix}\mid
BY+{}^tVA=0,\ {}^tYBY+AW+{}^tWA=0\right\},
$$
and
$$
U(F)=\left\{\begin{pmatrix}I&0&\begin{pmatrix}0&ex\\ -x&0\end{pmatrix}\\ 
0&I&0\\ 0&0&I\end{pmatrix}\mid x\in\RR\right\}. 
$$
\end{lemma}
\begin{proof}
  This is a straightforward calculation. 
\end{proof}
As in~\cite{Ko1} we realise $\cD_L$ as a Siegel domain and
$\cD_L(F)=U(F)_\CC\cD_L(F)$ is identified with
$\CC\times\CC^{n-2}\times \HH$. The identification is by choosing
homogeneous coordinates $(t_1\colon\ldots\colon t_{n+2})$ on
$\PP(L_\CC)$ so that $t_{n+2}=1$ and mapping $t_1\mapsto z\in \CC$,
$t_{n+1}\mapsto\tau\in\HH$ and $t_i\mapsto w_{i-2}\in\CC$ for $3\le i
\le n$: the value of $t_2$ is determined by the equation
\begin{equation}\label{siegeldomain}
  2\delta et_2=-2\delta z\tau-{}^t\underline{w}B\underline{w}
\end{equation}
where $\underline{w}\in\CC^{n-2}$ is a column vector. 

We are interested in the action of $N(F)_\ZZ=N(F)\cap \Gamma$ on
$\cD_L(F)$. We denote by ${\underline V}_i$ the $i$th row of the
matrix~$V$ in Lemma~\ref{describe_NUW_for_dim_1}. 
\begin{proposition}\label{action_of_NFZ}
If $g\in N(F)$ is given by
$$
\begin{pmatrix}U&V&W\\ 0&X&Y\\ 0&0&Z\end{pmatrix},\qquad
Z=\begin{pmatrix}a&b\\ c&d\end{pmatrix}
$$
then $g$ acts on $\cD_L(F)$ by
\begin{eqnarray*}
  z&\Mapsto& \frac{z}{\det Z}+(c\tau+d)^{-1}\bigg(\frac{c}{2\delta\det Z}
  {}^t\underline{w}B\underline{w}
  +\underline{V}_1\underline{w}+W_{11}\tau+W_{12}\bigg)\\
  \underline{w}&\Mapsto&(c\tau+d)^{-
    1}\left(X\underline{w}+Y\begin{pmatrix}\tau\\ 1\end{pmatrix}\right)\\
  \tau&\Mapsto& (a\tau+b)/(c\tau+d). 
\end{eqnarray*} 
\end{proposition}

\begin{proof}
  This is also a straightforward calculation. One need only take into
  account that
$$
U=\frac{1}{\det Z} \begin{pmatrix}d & -ce\\ b/e & a\end{pmatrix}
$$
\end{proof}

We must now describe $N(F)_\ZZ$ and $U(F)_\ZZ$. 
\begin{proposition}\label{Z_is_integral}
  If $g\in N(F)_\ZZ$ then $Z\in \SL_2(\ZZ)$, and if $g\in U(F)_\ZZ$
  then $x\in\ZZ$. 
\end{proposition}
\begin{proof}
  For $Z$, it is enough to show that $Z\in\Mat_{2\times 2}(\ZZ)$,
  since it acts on~$\HH$. The condition that $g\in N(F)_\ZZ$ or $g\in
  U(F)_\ZZ$ is that $N^{-1}gN\in \Gamma$ and in particular
  ${N'}^{-1}gN\in \GL_{n+2}(\ZZ)$. We calculate this directly:
$$
N^{-1}gN=\begin{pmatrix}U&V&VB^{-1}C+W-UT+TZ\\ 0&X&Y+XB^{-1}C
-B^{-1}CZ\\ 0&0&Z\end{pmatrix},
$$
so $Z$ is integral. In fact, because of ${}^tUAZ=A$ we even have $Z\in
\Gamma_0(e)$. 

If $g\in U(F)_\CC$ we have in addition $V=0$, $Y=0$, $U=Z=I_2$ and
$X=I_{n-2}$, so $VB^{-1}C+W-UT+TZ=W$ and therefore $W$ is integral. 
\end{proof}

Now we can calculate the action on the tangent space at a point in the
boundary. Suppose $g\in G(F)=N(F)_\ZZ/U(F)_\ZZ$ has finite order
$m>1$. We abuse notation by also using $g$ to denote a corresponding
element of~$N(F)_\ZZ$. We choose a coordinate $u=\exp_e(z):=e^{2\pi i
  z/e}$ on $U(F)_\CC/U(F)_\ZZ\imic \CC^*$, where $e$ is as in
Lemma~\ref{special_basis_for_L}, because $g\in U(F)_\ZZ$ acts by
$z\mapsto z+ex$. The compactification is given by allowing $u=0$. We
suppose that $g$ fixes the point $(0, \underline{w}_0,\tau_0)$. We
define $\Sigma(g)$ as we did before, in equation~(\ref{RST}), as
$\sum\{\frac{a_i}{m}\}$ if the eigenvalues are $\zeta^{a_i}$ for
$\zeta=e^{2\pi i/m}$. 

\begin{proposition}\label{RST_for_nonqrefs_at_dim1}
  If $n\ge 8$ and no power of $g$ acts as a quasi-reflection at $(0,
  \underline{w}_0,\tau_0)$ then $\Sigma(g)\ge 1$. 
\end{proposition}
\begin{proof}
  This closely follows \cite[(8.2)]{Ko1}. The action of $g$ on the
  tangent space is given by
$$
\begin{pmatrix}\exp_e(t)&0&0\\ *&(c\tau_0+d)^{-1}X&0\\ *&*&(c\tau_0+d)^{-
    2}\end{pmatrix}
$$ 
where $t=(c\tau_0+d)^{-1}(c{}^t\underline{w}_0 B\underline{w}_0/2
+\underline{V}_1\underline{w}_0+W_{11}\tau_0+W_{12})/e$, by
Lemma~\ref{action_of_NFZ}. Observe that $c\tau_0+d=\xi$ is a (not
necessarily primitive) fourth or sixth root of unity, because of the
well-known fixed points of $\SL_2(\ZZ)$ on~$\HH$. 

Suppose $X$ is of order~$m_X$. We consider the decomposition of the
representation $X$, {i.e.} of $E_\QQ^\perp/E_\QQ$ as a $g$-module. It
decomposes as a direct sum of $\cV_d$. If $\xi\neq\pm 1$ the situation
is exactly as in the case $\varphi(r)=2$ at the end of the proof of
Theorem~\ref{RST_for_nonqrefs_on_interior}, except that the right-hand
side of equation~(\ref{dimcount_phi=2}) is now equal to $n-2$ (that
is, $\rank X$) instead of~$n$. Any $\cV_d$ contributes at least
$c_{\min}(d)$ to $\Sigma(g)$, so we may assume that $\varphi(d)\le 2$;
but then the $1$-dimensional summands contribute at least
$\frac{1}{6}$ and the $2$-dimensional ones at least $\frac{1}{3}$. 
Moreover, if $m_X>2$ then $X$ has a pair of conjugate eigenvalues and
in the case $\xi=\pm 1$ they contribute $1$ to $\Sigma(g)$. 

So we may assume that $m_X=1$ or $m_X=2$, and $\xi=\pm 1$. Since
$-1\in\Gamma$ acts trivially on $\cD_L$ we may replace $g$ by $-g$ if
we prefer, and assume that $\xi=1$. Since $g$ fixes
$(0,\underline{w}_0,\tau_0)$ that implies $Z=I$. If also $m_X=1$, so
$X=I$, then by Proposition~\ref{action_of_NFZ} we have
$$
Y\begin{pmatrix}\tau_0\\ 1\end{pmatrix}=\underline{0}
$$
and since $\tau_0\not\in\ZZ$ this implies $Y=0$. But then ${}^tVA=0$
by Lemma~\ref{describe_NUW_for_dim_1}, so $g\in U(F)_\ZZ$. 

So the remaining possibility is that $Z=I$ and $m_X=2$: thus $U=I$ since
${}^tUAZ=A$, and $c=0$. But then $t$ is a half-integer, because
$$
\underline{w}_0=X\underline{w}_0+Y\begin{pmatrix}\tau_0\\ 1\end{pmatrix}
$$
and the condition $g^2\in U(F)_\ZZ$ implies that $VX=-V$, that $XY=-Y$
and that
$$
2W\equiv -VY \bmod \begin{pmatrix}0&e\\ -1&0\end{pmatrix}. 
$$
So, modulo~$e\ZZ$, we have 
\begin{eqnarray*}
2t&=&2\underline{V}_1\underline{w}_0+2W_{11}\tau_0+2W_{12}\\
&\equiv& 2\underline{V}_1\underline{w}_0-\underline{V}_1Y
\begin{pmatrix}\tau_0\\ 1\end{pmatrix}\\
&\equiv&\underline{V}_1(I+X)\underline{w}_0\\
&\equiv&0. 
\end{eqnarray*}
Thus the eigenvalue $\exp_e(t)$ is $\pm 1$, so in this case all
eigenvalues on the tangent space are $\pm 1$ and either $\Sigma(g)\ge
1$ or $g$ acts as a reflection. In particular any quasi-reflections
have order~$2$. 
\end{proof}

\begin{corollary}\label{no_fixed_divisors_at_dim1}
  There are no divisors at the boundary over a dimension~$1$ cusp~$F$
  that are fixed by a nontrivial element of~$G(F)$. 
\end{corollary}
\begin{proof}
  From the proof of Proposition~\ref{RST_for_nonqrefs_at_dim1}, any
  quasi-reflection $g$ has $m_X=2$, and hence fixes a divisor
  different from~$u=0$. 
\end{proof}
Finally we check the analogue of
Proposition~\ref{g^l_satisfies_modified_RST}. We define $\Sigma'(g)$
for $g\in G(F)$ exactly as in equation~(\ref{Sigmadash}). 
\begin{proposition}\label{g^l_satisfies_modified_RST_at_dim1}
  If $g\in G(F)$ is such that $g^k=h$ is a reflection and $n\ge 9$
  then $\Sigma'(g^l)\ge 1$ for every $1\le l<k$. 
\end{proposition}
\begin{proof}
  If the unique eigenvalue of $h$ that is different from~$1$ (hence
  equal to~$-1$) is $\exp_e(t)$ then the contribution from $X^l$ to
  $\Sigma'(g)$ is at least~$1$. Other\-wise, consider the $\cV_d$ (in
  the decomposition as a $g$-module) in which the exceptional
  eigenvector $\bde_0$ occurs, satisfying $h(\bde_0)=-\bde_0$. We must
  have $d=1$ or $d=2$, since if $\varphi(d)>1$ the eigenvalue $-1$ for
  $h$ would occur more than once. But the rest of $X$ ({i.e.} the
  $(n-3)$-dimensional $g$-module $E_\QQ^\perp/(E+\QQ\,\bde_0)$)
  contributes at least~$1$ to $\Sigma(g)$ and hence to $\Sigma'(g)$,
  as long as $n-3\ge 6$, as was shown in
  Proposition~\ref{RST_for_nonqrefs_at_dim1}. 
\end{proof}
\begin{corollary}\label{can_sings_on_dim1_cusps}
  If $n\ge 9$, the toroidal compactification $\Bar\cF_L(\Gamma)$ may
  be chosen so that on a boundary component over a dimension~$1$ cusp,
  $\Bar\cF_L(\Gamma)$ has canonical singularities, and there are no
  fixed divisors in the boundary. 
\end{corollary}
\begin{proof}
  This is immediate from Corollary~\ref{no_fixed_divisors_at_dim1},
  Proposition~\ref{RST_for_nonqrefs_at_dim1} and
  Proposition~\ref{g^l_satisfies_modified_RST_at_dim1}. In fact there
  are no choices to be made in this part of the boundary. 
\end{proof}

\section{Special reflections in $\Tilde{\Orth}(L)$}\label{reflections}

Let $L$ be an arbitrary nondegenerate integral lattice, and
write $D$ for the exponent of the finite group $A_L=L^\vee/L$. 
The reflection with respect to the hyperplane defined
by a vector $r$ is given by
$$
\sigma_r\colon l\Mapsto l-\frac{2(l,r)}{(r,r)}r. 
$$
For any $l\in L$ its {\it divisor\/} $\divv(l)$ in $L$ is the positive
generator of the ideal $(l,L)$.  In other words $l^*=l/\divv(l)$ is a
primitive element of the dual lattice $L^\vee$.  If $r$ is primitive
and the reflection $\sigma_r$ fixes $L$, {i.e.} $\sigma_r\in
\Orth(L)$, then we say that $r$ is a reflective vector. In this case
\begin{equation}\label{refldiv}
\divv (r)\mid r^2 \mid 2\divv (r). 
\end{equation}
\begin{proposition}\label{reflid}
  Let $L$ be a nondegenerate even integral lattice.  Let $r\in L$ be
  primitive.  Then $\sigma_r\in \widetilde \Orth(L)$ if and only if
  $r^2=\pm 2$. 
\end{proposition}
\begin{proof} For $r^*=r/\divv(r)\in L^\vee$ and $\sigma_r\in
  \Tilde\Orth(L)$ we get
$$
\sigma_r(r^*)=-r^*\equiv r^*\mod L. 
$$
Therefore $2r^*\in L$, $\divv(r)=1$ or $2$ (because $r$ is primitive)
and $r^2=\pm 2$ or $\pm 4$, because $L$ is even.  If $r^2=\pm 2$ then
$\sigma_r\in \widetilde \Orth(L)$.  If $r^2=\pm 4$, then $\divv(r)=2$
by condition~(\ref{refldiv}).  For such $r$ the reflection $\sigma_r$
is in $\widetilde \Orth(L)$ if and only if
$$
l^\vee-\sigma_r(l^\vee)=\frac{(r,l^\vee)}{2}r=(r^*,l^\vee)r\in L
$$
for any $l^\vee\in L^\vee$. Therefore $r^*=r/2\in (L^\vee)^\vee=L$. 
We obtain a contradiction because $r$ is primitive. 
\end{proof}

\begin{proposition}\label{reflminusid}
  Let $L$ be as in Proposition~\ref{reflid} and let $r\in L$ be
  primitive.  If $-\sigma_r\in \widetilde \Orth(L)$, {i.e.} 
  $\sigma_r|_{A_L}=-\id$, then
\begin{itemize}
\item[(i)] $r^2=\pm 2D$ and $\divv(r)=D\equiv 1\mod 2$, or $r^2=\pm D$
  and $\divv(r)=D$ or $D/2$;
\item[(ii)] $A_L\cong (\ZZ/2\ZZ)^m\times (\ZZ/D\ZZ)$. 
\end{itemize}
In the opposite direction we have
\begin{itemize}
\item[(iii)] If $r^2=\pm D$ and either $\divv(r)=D$ or
  $\divv(r)=D/2\equiv 1\mod 2$, then $-\sigma_r\in \widetilde
  \Orth(L)$;
\item[(iv)] If $r^2=\pm 2D$ and $\divv(r)=D\equiv 1\mod 2$, then
  $-\sigma_r\in \widetilde \Orth(L)$. 
\end{itemize}
\end{proposition}
\begin{proof}
  (i) $\sigma_r|_{A_L}=-\id$ is equivalent to the following condition:
\begin{equation}\label{-id}
  2l^\vee\equiv \frac{2(r,l^\vee)}{r^2}r\mod L
  \qquad\forall\ l^\vee\in L^\vee. 
\end{equation}
It follows that if $r^2=2e$, then $(2L^{\vee})/ L$ is a subgroup of
the cyclic group $\langle(r/e)+L\rangle$.  Thus $D$ divides $2e$.  But
by definition of the divisor of the vector $e\mid \divv(r)\mid D$,
therefore
$$
e\mid \divv(r)\mid 2e\quad\text{and}\quad e\mid D\mid 2e. 
$$
From this it follows that $(2L^\vee)/L$ is a subgroup of the cyclic
group generated by $(r/D)+L$ or $(2r/D)+L$. This implies~(ii). 

Let us assume that $r^2=\pm 2D$ and $\divv(r)=D\equiv 0\mod 2$.  We
have $2l^\vee\equiv \pm \frac{(r,l^\vee)}{D}r\mod L$.  If the order of
$l^\vee$ in the discriminant group is odd, then $(r,l^\vee)$ is even,
since $D$ is even.  If the order of $l^\vee$ is even, then
$(r,l^\vee)$ is again even, because the order of $2l^\vee$ is $D/2$. 
Therefore $(r/2,l^\vee)\in \ZZ$ for all $l^\vee\in L^\vee$.  This
contradicts the assumption that $r$ is primitive.  Thus (i) is proved. 

(iii) Let assume that $\divv(r)=D$.  In this case $r^*=r/D$ and $
2r^*+L$ is a generator of $(2L^\vee)/L$. According to (ii) we have
that for any $l^\vee\in L^\vee$, $2l^\vee= 2xr^*+l'$, where $x\in
\ZZ$, $l'\in L$.  Therefore
\begin{equation}\label{DD}
  \frac{(2l^\vee,r)}{r^2}r=
  2xr^*\pm \frac{(l',r)}{D}r
  \equiv
  2xr^*\equiv 2l^\vee\mod L
\end{equation}
and $-\sigma_r\in \widetilde \Orth(L)$ according to
condition~(\ref{-id}). 

Let assume that $\divv(r)=D/2\equiv 1 \mod 2$.  We have to check
condition~(\ref{-id}) for all elements of order $2$ or $D$ in $A_L$. 
If $\ord(l^\vee)=2$, then $(2l^\vee,r)\equiv 0\mod D/2$, and also
$(l^\vee,r)\equiv 0\mod D/2$, because $D/2$ is odd.  It follows that
$2(l^\vee,r)/r^2\in \ZZ$.  If $l^\vee$ is an element of order $D$, we
have $2l^\vee= 2xr^*+l'$ as above with $r^*=(2r)/D$ and $l'\in L$. 
Thus $(l',r)$ is even. But $(l',r)$ is also divisible by the odd
number $D/2$.  Therefore $(l',r)\equiv 0\mod D$ and
equation~(\ref{DD}) is also true. 

(iv) is similar to (iii). $D$ is odd and the group $A_L$ is cyclic
with generator $r^*=r/D$.  Therefore $l^\vee=xr^*+l'$ for any
$l^\vee\in L^\vee$ and
$$
\frac{(2l^\vee,r)}{r^2}r=\frac{2(xr^*+l',r)}{r^2}r=
2xr^*\pm \frac{2(l',r)}{2D}\equiv 2l^\vee\mod L. 
$$
\end{proof}
\begin{corollary}\label{odddet}
  Let $L$ be an even integral lattice and $|A_L|=|\det L|$ be odd. 
  Then
\begin{itemize}
\item[(i)] $\sigma_r\in \widetilde \Orth(L)$ if and only if $r^2=\pm
  2$;
\item[(ii)] $-\sigma_r\in \widetilde\Orth(L)$ if and only if $r^2=\pm
  2D$ and $\divv(r)=D$. 
\end{itemize}
\end{corollary}
With $\Kthree$ surfaces in mind, we consider in more detail the
lattice $L_{2d}=2U\oplus 2E_8(-1)\oplus \latt{-2d}$. 
\begin{corollary}\label{reflK3}
  Let $\sigma_r$ be a reflection in $\Orth(L_{2d})$ defined by a
  primitive vector $r\in L_{2d}$.  $\sigma_r$ induces $\pm\id$ on the
  discriminant form $L_{2d}^\vee/L_{2d}$ if and only if $r^2=\pm 2$ or
  $r^2=\pm 2d$ and $\divv(r)=d$ or $2d$. 
\end{corollary}
\begin{proof}
  Any $r\in L_{2d}$ can be written as $r=m+xh$, where $m\in
  L_0=2U\oplus2E_8(-1)$ and $h^2=-2d$ ($h$ is primitive). 

If $r^2=\pm 2d$ and $\divv(r)=2d$, then $-\sigma_r\in
\Tilde{\Orth}(L_{2d})$ by Proposition~\ref{reflminusid}. 

If $r^2=\pm 2d$ and $\divv(r)=d$, then $r=dm_0+xh$, where
$x^2=1-d(m_0^2/2)$.  We see that
$$
\sigma_r\left(\frac{h}{2d}\right)=\frac{h}{2d}(1-2x^2)-xm_0\equiv
 -\frac{h}{2d}\mod L_{2d}. 
$$
\end{proof}

The types of reflections in the full orthogonal group $\Orth^+(L)$ for
$L=L_{2d}^{(0)}=2U\oplus \latt{-2d}$ were classified in \cite{GH2}
(for square-free~$d$).  The result for $L_{2d}=2U\oplus 2E_8(-1)\oplus
\latt{-2d}$ is exactly the same, because the unimodular part
$2E_8(-1)$ plays no role in the classification.

The reflection $\sigma_r$ is an element of $\Orth^+(L_\RR)$ (where $L$
has signature $(2,n)$) if and only if $r^2<0$: see \cite{GHS1}. 

The $(-2)$-vectors of $L_{2d}$ form one or two (if $d\equiv 1\mod 4$)
orbits with respect to the group $\widetilde \Orth^+(L_{2d})$.  We can
also compute the number of $\widetilde\Orth^+(L_{2d})$-orbits of the
$(-2d)$-reflective vectors in Corollary~\ref{reflK3}.  However, in
this paper we only need to know the orthogonal complements of
$(-2d)$-vectors, which we compute in Proposition~\ref{-2dvect}.  (For
the case of $(-2)$-vectors see \cite[\S{}3.6]{GHS1}). 

The following lemma, which we use in the proof of
Proposition~\ref{-2dvect}, is well-known, but we state it and give a
general proof here for the convenience of the reader.  Recall that an
integral lattice $T$ is called $2$-elementary if $A_T=T^\vee/T\cong
(\ZZ/ 2\ZZ )^m$. 

\begin{lemma}\label{2-elementary}
  Let $T$ be a primitive sublattice of an unimodular even lattice $M$,
  and let $S$ be the orthogonal complement of $T$ in $M$. Suppose that
  there is an involution $\sigma\in \Orth(M)$ such that
  $\sigma|_T=\id_T$ and $\sigma|_S=-\id_S$. Then $T$ and $S$ are
  $2$-elementary lattices. 
\end{lemma}
\begin{proof}
  Let us consider the inclusions $T\oplus S\subset M\subset T^\vee
  \oplus S^\vee$.  We have that $(A_S,\ q_{S})\cong (A_T,\ -q_{T})$
  because $M$ is unimodular (see \cite{Nik2}).  In particular
  $[M:T\oplus S]=[T^\vee:T]=[S^\vee:S]$.  It follows that
$$
H=M/(T\oplus S)\cong \phi(M)/S=S^\vee/S=A_S. 
$$
Here $\phi\colon M\to S^\vee$ is defined by $\phi(m)(s)=(m,s)$ where
$s\in S$.  The natural projections of the subgroup $H<A_T\oplus A_S$
onto $A_T$ and $A_S$ are injective, therefore the action of $\sigma$
on $A_S$ is completely determined by the action of $\sigma$ on $A_T$. 
Thus $\sigma$ acts trivially on $A_S$ since it acts trivially on
$A_T$.  But we assumed that $\sigma(s^\vee)=-s^\vee$ for any
$s^\vee\in S^\vee$.  It follows that $A_S$ is an abelian $2$-group. 
\end{proof}

\begin{proposition}\label{-2dvect}
  Let $r$ be a primitive vector of $L_{2d}$. If $\divv(r)=2d$ then
$$
r^\perp_{L_{2d}}\cong 2U\oplus 2E_8(-1). 
$$
If $\divv(r)=d$ then either
$$
r^\perp_{L_{2d}}\cong U\oplus 2E_8(-1)\oplus \latt{2} \oplus \latt{-2}
$$
or
$$
r^\perp_{L_{2d}}\cong U\oplus 2E_8(-1)\oplus U(2). 
$$
\end{proposition}
\begin{proof}
  The lattice $L_{2d}$ is the orthogonal complement of a primitive
  vector $h$, with $h^2=2d$ in the unimodular $\Kthree$ lattice
  $L_{\Kthree}=3U\oplus 2E_8(-1)$.  We put $L_r=r^\perp_{L_{2d}}$ and
  $S_r=(L_r)^\perp_{L_{\Kthree}}$. 

  We note that $L_r$ and $S_r$ have the same determinant: in fact
$$
\det L_r= \det S_r=4d^2/\divv(r)^2=
\begin{cases}
1&\quad\text{if } \divv(r)=2d,\\
4&\quad\text{if } \divv(r)=d. 
\end{cases}
$$ 
To see this, consider a more general situation.  Let $N$ be a
primitive even nondegenerate sublattice of any even integral lattice
$L$ and let $N^\perp$ be its orthogonal complement in $L$.  Then we
have
$$
N\oplus N^\perp\subset L\subset L^\vee
\subset N^\vee\oplus (N^\perp)^\vee,
$$
where $L/(N\oplus N^\perp)\cong L^\vee/(N^\vee\oplus (N^\perp)^\vee)$. 
As before we have $\phi\colon L\to N^\vee$, and $\ker(\phi)=N^\perp$. 
Since $L/(N\oplus N^\perp)\cong \phi(L)/N$ we obtain
$$
|L/(N\oplus N^\perp)|= |\phi(L)/N|=|\det N|/[N^\vee:\phi(L)],
$$
as $|\det N|=[N^\vee:N]$.  From the inclusions above
$$
|\det N|\cdot|\det N^\perp|=(|\det L|)[\phi(M):N]^2=|\det L|\cdot|\det
N|^2/[N^\vee:\phi(L)]^2. 
$$

In our particular case $L=L_{2d}$, $N=\ZZ r$ and $L_r=N^\perp$.  We
have $[N^\vee:\phi(L)]=\divv(r)$, where $\divv(r)\ZZ=(r,L)$, and this
gives us the formula for the determinant of $L_r$. 

If $\divv(r)=2d$ then $L_r$ and $S_r$ are are isomorphic to the unique
unimodular lattices of signatures $(2,18)$ and $(1,1)$ respectively:
that is, $L_r\cong 2U\oplus 2E_8(-1)$ and $S_r\cong U$. 

If $\divv(r)=d$ then the reflection $\sigma_r$ acts as $-\id$ on the
discriminant group (see Corollary~\ref{reflK3}).  Therefore we can
extend $-\sigma_r\in \widetilde \Orth(L_{2d})$ to an element of
$\Orth(L_{\Kthree})$ by putting $(-\sigma_r)|_{\ZZ h}=\id$.  So
$\sigma_r$ has an extension $\tilde\sigma_r\in \Orth(L_{\Kthree})$
such that $\tilde\sigma_r|_{L_r}=\id_{L_r}$ and
$\tilde\sigma_r|_{S_r}=-\id_{S_r}$.  It follows from
Lemma~\ref{2-elementary} that $L_r$ and $S_r$ are $2$-elementary
lattices. 

The finite discriminant forms of $2$-elementary lattices were
classified by Nikulin in~\cite{Nik3}.  The genus of $M$ (and the class
of $M$ if $M$ is indefinite) is determined by the signature of $M$,
the number of generators $m$ of $A_M$ and the parity $\delta_M$ of the
finite quadratic form $q_M\colon A_M\to \QQ/2\ZZ$, which is given by
$\delta_M= 0$ if $l^2\in \ZZ$ for all $l\in M^\vee$ and $\delta_M=1$
otherwise: (see [Nik3, \S{}3]).  In particular, for an indefinite
lattice $S_r$ of rank~$2$ and determinant~$4$ we have
$$
S_r\cong
\begin{cases}\ U(2)&\quad\text{if}\quad \delta_{S_r}=0,\\
  \ \latt 2 \oplus \latt{-2}& \quad\text{if}\quad \delta_{S_r}=1. 
\end{cases}
$$
The class of the indefinite lattice $L_r$ is uniquely defined by its
discriminant form.  Proposition~\ref{-2dvect} is proved. 
\end{proof}

Geometrically the three cases in Proposition~\ref{-2dvect} correspond
to the N\'eron-Severi group being (generically) $U$, $U(2)$ or
$\latt{2}\oplus\latt{-2}$ respectively. The $\Kthree$ surfaces
(without polarisation) themselves are, respectively, a double cover of
the Hirzebruch surface $F_4$, a double cover of a quadric, and the
desingularisation of a double cover of $\PP^2$ branched along a nodal
sextic.

\section{Special cusp forms.}\label{specialcusp}

Let $L=2U\oplus L_0$ be an even lattice of signature $(2,n)$ ($n\ge
3$) containing two hyperbolic planes.  We write
$\cF_L=\cF_L(\widetilde\Orth^+(L))$ for brevity.  A $0$-dimensional
cusp of $\cF_L$ is defined by a primitive isotropic vector $v$.  Any
two primitive isotropic vectors of divisor $1$ lie in the same
$\Tilde\Orth^+(L)$-orbit, according to the well-known criterion of
Eichler (see \cite[\S 10]{E}).  We call the corresponding cusp the
{\it standard $0$-dimensional cusp\/} of the Baily--Borel
compactification $\cF_L^{*}$. 

Each $1$-dimensional boundary component $F$ of $\cD_L$ is isomorphic
to the upper half plane $\HH$ and in the Baily--Borel compactification
this corresponds to adding an (open) curve $\Lambda \backslash \HH$,
where $\Lambda\subset \SL_2(\QQ)$ is an arithmetic group which depends
on the component $F$. Details of this can be found in \cite{BB} and
\cite{Sc}. For our purpose we need one general result not contained
there. 

\begin{lemma}\label{cuspclosure} 
  Suppose that $L$ is even, and that any isotropic subgroup of the
  discriminant group $(A_L,q_L)$ is cyclic.  Then the closure of every
  $1$-dimensional cusp in $\cF_L^*$ contains the standard
  $0$-dimensional cusp. 
\end{lemma}
\begin{proof}
  Let $E$ be a primitive totally isotropic rank $2$ sublattice of $L$
  and define the lattice $\widetilde E=E_{L^\vee}^{\perp\perp}$ (both
  orthogonal complements are taken in the dual lattice $L^{\vee}$). 
  We remark that $E\subset \widetilde E$ and that $E=\widetilde E\cap
  L$ because $E$ is isotropic and primitive.  Thus the finite group
$$
H_E= E_{L^\vee}^{\perp\perp}/E < A_L
$$
is an isotropic subgroup of the discriminant group of $L$.  Let us
take a basis of $L$ as in Lemma~\ref{special_basis_for_L}.  It is easy
to see that
$$
H_E \cong A^{-1} \ZZ^2/\ZZ^2. 
$$ 
In the case we are considering, $H_E$ is a cyclic subgroup ($|H_E|^2$
divides $\det L$).  Therefore $A=\diag(1,e)$.  Thus $E$ contains
primitive isotropic vectors with divisors $1$ and $e$, and the first
vector defines the standard $0$-dimensional cusp. 
\end{proof}

\noindent{\it Remark.} 
If the discriminant group of $L$ contains a non-cyclic isotropic
subgroup then there is a totally isotropic sublattice $E$ of $L$ such
that the finite abelian group $H_E$ has elementary divisors
$(\delta,\delta e)$ with $\delta>1$.  Thus $\det L$ is divisible by
$\delta^4e^2$. 

Let $L=2U\oplus L_0$ be of signature $(2,n)$ and $u$ be a primitive
isotropic vector of divisor $1$. The tube realisation $\cH_u$ of the
homogeneous domain $\cD_L$ at the standard $0$-dimensional cusp is
defined by the sublattice $L_1=u^{\perp}/\ZZ u\cong U\oplus L_0$:
\begin{equation}\label{tube}
  \cH_u=\cH(L_1)=\{Z\in L_{1}\otimes \CC\ |\
  (\Im Z, \Im Z)>0\}^+,
\end{equation}
where ${}^+$ denotes a connected component of the domain (see \cite{G}
for details).  The modular group $\Tilde\Orth^+(L)$ acting on
$\cH(L_1)$ contains all translations $Z\to Z+l$ ($l\in L_1$). 
Therefore the Fourier expansion of a $\Tilde\Orth^+(L)$-modular form
$F$ at the standard cusp is
\begin{equation}\label{fourier}
  F(Z)=
  \sum_{l\in L_1^\vee,\ (l,l)\ge 0}
  a(l)\exp(2\pi i (l, Z)). 
\end{equation}

\begin{theorem}\label{cusp}
  Let $L$ be an even lattice with two hyperbolic planes such that any
  isotropic subgroup of the discriminant group of $L$ is cyclic. Let
  $F$ be a modular form with respect to $\Tilde\Orth^+(L)$. If its Fourier
  coefficients $a(l)$ at the standard cusp satisfy $a(l)=0$ if
  $(l,l)=0$, then $F$ is a cusp form. 
\end{theorem}
\begin{proof}
  A standard $1$-dimensional cusp is defined by a primitive totally
  isotropic sublattice $E_1=\latt{u, v}$ with $\divv(u)=\divv(v)=1$. 
  We can choose $(u,v)$ in such a way that they generate the maximal
  totally isotropic sublattice in $U\oplus U$.  Let $E$ be an
  arbitrary primitive totally isotropic sublattice of rank~$2$ of $L$
  defining a $1$-dimensional cusp of $\cF_L$.  We can assume that
  $E=\latt{u, v'}_\ZZ$ where $u$ defines the standard $0$-dimensional
  cusp (see Lemma~\ref{cuspclosure} above).  According to the Witt
  theorem for the rational hyperbolic quadratic space $L_1\otimes \QQ$
  there exists $\sigma \in \Orth(L_1\otimes \QQ)$ such that
  $\sigma(v')=v$.  We can extend $\sigma$ to an element of
  $\Orth^+(L\otimes \QQ)$ by putting $\sigma(u)=\pm u$.  The Siegel
  operator $\Phi_E$ for the boundary component defined by $E$ has the
  property $\Phi_E(F\circ \sigma)=\Phi_{\sigma(E)}(F)\circ \sigma$
  (see \cite{BB}).  Therefore
$$
\Phi_E(F)
=\Phi_{\sigma^{-1}E_1}(F)
=\Phi_{E_1}(F\circ \sigma^{-1})\circ \sigma. 
$$
We can calculate the Fourier expansion of the function under the
Siegel operator $\Phi_{E_1}$:
\begin{eqnarray}
  F\circ \sigma^{-1}&=&
  \pm \sum_{l\in L_1^\vee,\,(l,l)>0}
  a(l)\exp(2\pi i(l, \sigma^{-1}Z))\nonumber\\
  &=&\pm\sum_{l_1\in \sigma L_1^{\vee},\,(l_1,l_1)>0}
  a(\sigma^{-1}l_1)\exp(2\pi i(l_1, Z)). 
\end{eqnarray}
Thus $\Phi_{E}(F)= \Phi_{E_1}(F\circ \sigma^{-1})\circ \sigma\equiv 0$
and $F$ is a cusp form. 
\end{proof}

In \cite[Theorem 3.1]{G} modular forms for $\widetilde{\SO}^+(L)$ are
constructed using the arithmetic lifting of a Jacobi form $\phi$.  The
modular form $\Lift(\phi)$ is defined by its first Fourier-Jacobi
coefficient at a fixed standard $1$-dimensional cusp.  In particular,
we know the Fourier expansion at the standard $0$-dimensional cusp. 
Therefore we obtain the following improvement of the result proved in
\cite{G} for square-free $d$. 
\begin{corollary}\label{SOF} 
  Let $L=L_{2d}=2U\oplus 2E_8(-1)\oplus \latt{-2d}$.  Then the
  arithmetic lifting $\Lift(\phi)$ of a Jacobi cusp form $\phi\in
  J_{k,1}^{\cusp}(L_{2d})$ of weight $k$ and index~$1$ is a cusp form
  of weight $k$ for $\widetilde{\SO}^+(L_{2d})$ for any $d\ge 1$. 
\end{corollary}

\section{Application: K3 surfaces with a spin structure}\label{spinK3}

Instead of $\widetilde{\Orth}^+(L_{2d})$ and $\cF_{2d}$, we may
consider the subgroup $\widetilde{\SO}^+(L_{2d})$ of
$\widetilde{\Orth}^+(L_{2d})$ of index~$2$ and the corresponding
quotient
$$
{\cS\cF}_{2d} = \widetilde{\SO}^+(L_{2d}) \backslash \cD_{L_{2d}}. 
$$ 
If $d>1$ then ${\cS\cF}_{2d}$ is a double covering of ${\cF}_{2d}$. 
(For $d=1$ the two spaces coincide since
$\widetilde{\SO}^+(L_{2})\cong \widetilde{\Orth}^+(L_{2})/\pm I$.) 
This double covering has the following geometric interpretation: the
domain $\cD_{L_{2d}}$ is the parameter space of marked $\Kthree$
surfaces of degree $2d$, and dividing out by the group
$\widetilde{\Orth}^+(L_{2d})$ identifies all the different markings on
a given $\Kthree$ surface. Two markings will be identified under the
group $\widetilde{\SO}^+(L_{2d})$ if and only if they have the same
orientation. Hence ${\cS\cF}_{2d}$ parametrises polarised $\Kthree$
surfaces $(S,h)$ together with an orientation of the lattice
$L_h=h^{\perp}$. We shall refer to these as {\em oriented\/} $\Kthree$
surfaces. An orientation on a surface $S$ is also sometimes called a
{\em spin structure\/} on $S$. 

We have seen in Corollary~\ref{reflK3} that the branch divisor of the
map $\cD_{L_{2d}} \to \cF_{2d}$ is given by the divisors associated to
reflections $\sigma_r$ defined by a primitive vector $r$ of length
either $r^2=-2$ or $r^2=-2d$. Note that in the first case $\sigma_r$
acts trivially on the discriminant group whereas it acts as $-\id$ in
the second case. Hence $\pm \sigma_r \notin \widetilde{\SO}^+(L_{2d})$
if $r^2=-2$, but $-\sigma_r \in \widetilde{\SO}^+(L_{2d})$ if
$r^2=-2d$.  It follows that the quotient map $\cD_{L_{2d}} \to
{\cS\cF}_{2d}$ is branched along the $(-2d)$-divisors whereas the
double cover $ {\cS\cF}_{2d} \to {\cF}_{2d}$ is branched along the
$(-2)$-divisors. In this way the group $\widetilde{\SO}^+(L_{2d})$
separates the two types of contributions to our reflective
obstructions.  The reflective obstructions coming from the $(-2d)$
divisors are less problematic, as we shall see in the next theorem. 
The $(-2d)$-divisors have a geometric interpretation. The general
point on such a divisor is associated to a $\Kthree$ surface $S$ whose
transcendental lattice $T_S$ has rank $20$ and which admits an
involution acting as $-\id$ on $T_S$. For $d=p^2$ this was shown in
(\cite[Prop. 7.4]{Ko1}), and for general $d$ it follows from
Corollary~\ref{reflK3} and the proof of Proposition~\ref{-2dvect}
above. 

In \cite{G} it was proved that the modular variety
$\widetilde{\SO}^+(L_{2d})(q)\backslash \cD_{L_{2d}}$, where
$\widetilde{\SO}^+(L_{2d}(q))$ is the principal congruence subgroup of
$\widetilde{\SO}^+(L_{2d})$ of level $q\ge 3$, is of general type for
any $d\ge 1$.  Here we obtain a much stronger result. 

\begin{theorem}\label{orientedK3}
  The moduli space ${\cS\cF}_{2d} = \widetilde{\SO}^+(L_{2d})
  \backslash \cD_{L_{2d}}$ of oriented $\Kthree$ surfaces of degree
  $2d$ is of general type if $d\ge 3$. 
\end{theorem}
\begin{proof}
  For $L_{2d}=2U\oplus 2E_8(-1)\oplus \latt{-2d}$ the corresponding
  space of Jacobi cusp forms in $18$ variables is isomorphic (as a
  linear space) to the space of Jacobi cusp forms of Eichler-Zagier
  type (see \cite[lemma 2.4]{G})
$$
J_{k,1}^{\cusp}(L_{2d})\cong J_{k-8,d}^{\cusp}(EZ). 
$$ 
For $k=17$, this space is non-trivial for any $d\ge 3$.  Therefore for
any $d\ge 3$ there is a cusp form $F_{17}$ of weight $17$ with respect
to $\widetilde{\SO}^+(L_{2d})$. 

The ramification divisor of the projection $\pi_{\SO}\colon
\cD_{L_{2d}}\to \widetilde{\SO}^+(L_{2d})\backslash \cD_{L_{2d}}$ is
defined by $(-2d)$-reflections of $L_{2d}$. In Lemma~\ref{oddk} below
we show that the cusp form $F_{17}$ vanishes on the ramification
divisors of $\pi_{\SO}$. 

Hence $\cS\cF_{2d}$ is of general type for $d\ge 3$ by
Theorem~\ref{general_gt}. 
\end{proof}

\begin{lemma}\label{oddk}
  Any modular form $F\in M_{2k+1}(\widetilde{\SO}^+(L_{2d}))$ of odd
  weight vanishes along the divisors defined by $(-2d)$-reflective
  vectors. 
\end{lemma}
\begin{proof}
  Let $\sigma_r\in \Orth^+(L_{2d})$ be a reflection with respect to a
  $(-2d)$-vector. Then $-\sigma_r\in \widetilde{\SO}^+(L_{2d})$ (see
  Corollary~\ref{reflK3}).  For any $z\in \cD_{L_{2d}}$ with $(z,r)=0$
  and a modular form $F\in M_{2k+1}(\widetilde{\SO}^+(L_{2d}))$ we
  have
$$
F(z)=F((-\sigma_{r})(z))=F(-z)=(-1)^{2k+1}F(z),
$$
so $F(z)\equiv 0$. 
\end{proof}

We note that ${\cS\cF}_{2}={\cF}_{2}$ is unirational. 

The geometric interpretation of the $(-2)$-divisors, which form the
ramification of the covering ${\cS\cF}_{2d} \to {\cF}_{2d}$, is that
they parametrise those polarised $\Kthree$ surfaces whose polarisation
is only semi-ample, but not ample.  This is due the presence of
rational curves on which the polarisation has degree~$0$.  Thus in the
case $d=2$ the map ${\cS\cF}_{4} \to {\cF}_{4}$ is a double cover of
the moduli space of quartic surfaces branched along the discriminant
divisor of singular quartics.  The variety ${\cF}_{4}$ is unirational
but ${\cS\cF}_{4}$ is not, since there exists a canonical differential
form on it (see \cite{G}). There is also a cusp form of weight~$18$
with respect to $\widetilde{\SO}^+(L_{4})$ which vanishes on one of
the two irreducible components of the ramification divisors for $d=2$. 
We shall return to this question in a more general context in
\cite{GHS2}.

\section{Pull-back of the Borcherds function $\Phi_{12}$.}\label{Borcherds}

To construct pluricanonical differential forms on a smooth model of
$\cF_{2d}$ we shall use the pull-back of the Borcherds automorphic
product $\Phi_{12}$. 

Let $L_{2,26}=2U\oplus 3E_8(-1)$ be the unimodular lattice of
signature $(2,26)$. For later use, we note the following simple lemma. 
\begin{lemma}\label{orthcomp2d}
  Let $r$ be a primitive reflective vector in $L_{2d}$ with $r^2=-2d$
  and let $L_r=r^\perp_{L_{2d}}$ be its orthogonal complement
  considered as a primitive sublattice of the unimodular lattice
  $L_{2,26}$.  Then
$$
(L_r)^\perp_{L_{2,26}}\cong E_8(-1), \ E_7(-1)\oplus\latt{-2} \ \text{
  or }\ D_8(-1). 
$$
\end{lemma}
\begin{proof} 
  In the proof of Proposition \ref{-2dvect} we found $L_r$ and its
  orthogonal complement $S_r$ in the unimodular lattice
  $L_{\Kthree}=3U+2E_8(-1)$.  The discriminant forms of $S_r$ and
  $K_r=(L_r)^\perp_{L_{2,26}}$ coincide, but $K_r$ is of signature
  $(0,8)$.  The three possible genera of $K_r$ are represented by
  $E_8(-1)$, $E_7(-1)\oplus \langle -2\rangle$ and $D_8(-1)$.  The
  genera of such lattices contain only one class: one can can prove
  this well-known fact by analysing sublattices of order $2$ in $E_8$
  or simply check it using MAGMA. 
\end{proof}

The Borcherds function $\Phi_{12}\in M_{12}(\Orth^+(L_{2,26}),\det)$
is the unique modular form of weight $12$ and character $\det$ with
respect to $\Orth^+(L_{2,26})$ (see \cite{B}). 

$\Phi_{12}$ is the denominator function of the fake Monster Lie
algebra and it has a lot of remarkable properties.  In particular, the
zeros of $\Phi_{12}(Z)$ lie on rational quadratic divisors defined by
$(-2)$-vectors in $L_{2,26}$, {i.e.}, $\Phi_{12}(Z)=0$ if and only if
there exists $r\in L_{2,26}$ with $r^2=-2$ such that $(r,Z)=0$ and the
multiplicity of the rational quadratic divisor in the divisor of zeros
of $\Phi_{12}$ is~$1$. 

Pulling back this function gives us many interesting automorphic forms
(see \cite[pp. 200-201]{B}, \cite[pp. 257-258]{GN}).  In the context
of the moduli of $\Kthree$ surfaces this function was used in
\cite{BKPS} and \cite{Ko2}.  We summarise their results in a suitable
form. 

Let $l\in E_8(-1)$ satisfy $l^2=-2d$. The choice of $l$ determines an
embedding of $L_{2d}$ into $L_{2,26}$ as well as an embedding of the
domain $\cD_{L_{2d}}$ into $\cD_{L_{2,26}}$.  We put $R_l=\{r\in
E_8(-1)\mid r^2=-2,\ (r, l)=0\}$, and $N_l=\# R_l$. (It is clear that
$N_l$ is even.)  Then by \cite{BKPS} the function
\begin{equation}\label{pb}
  \left. F_l=
    \frac{\Phi_{12}(Z)}{
      \prod_{\{ \pm r\}\in R_l} (Z, r)}
    \ \right\vert_{\cD_{L_{2d}}}
  \in M_{12+\frac{N_l}2}(\Tilde\Orth^+(L_{2d}),\, \det)
\end{equation}
is a non-trivial modular form of weight $12+\frac{N_l}2$ vanishing on
all $(-2)$-divisors of $\cD_{L_{2d}}$. (As we did in
Section~\ref{specialcusp}, we think of a modular form as a function on
$\cD_L$ rather than $\cD_L^\bullet$, by identifying $\cD_L$ with a
tube domain realisation as in equation~(\ref{tube}) above.)  Moreover
it is shown in~\cite{Ko2} that $F_l$ is a cusp form if $d$ is
square-free and the weight is odd. 

In fact much more is true. 
\begin{theorem}\label{cuspform} 
  The function $F_l$ has the following properties:
\begin{itemize}
\item[(i)] $F_l\in M_{12+\frac{N_l}2}(\Tilde{\Orth}^+(L_{2d}),\,
  \det)$ and $F_l$ vanishes on all $(-2)$-divisors. 
\item[(ii)] $F_l$ is a cusp form for any $d$ if $N_l>0$. 
\item[(iii)] If the weight of $F_l(Z)$ is smaller than $68$ ({i.e.},
  if $N_l<112$) then $F_l(Z)$ is zero along the branch divisor of the
  projection
$$
\pi\colon \cD_{L_{2d}}\To \Gamma_{2d}\backslash \cD_{L_{2d}}=\cF_{2d}. 
$$
\end{itemize}
\end{theorem}
\begin{proof} 
  (i) was proved in \cite{BKPS}, but we repeat some details here for
  convenience. First, $F_l(Z)$ is holomorphic because of the properties
  of the divisor of $\Phi_{12}$. Then $F_l(tZ)=t^{-(12+N_l/2)}F_l(Z)$
  for any $Z\in \cD_{L_{2d}}$.  Any $g\in \widetilde \Orth^+(L_{2d})$
  can be extended (by the identity on the orthogonal complement of
  $L_{2d}$ in $L_{2,26}$) to an element $\tilde g$ of
  $\Orth^+(L_{2,26})$.  Therefore $F_l(gZ)=\det(g)F_l(Z)$ since
  $\tilde g(r)=r$ for all roots in $R_l$. This modular form is
  evidently not identically zero. On the other hand, because it has
  character $\det$ it vanishes on all divisors of $\cD_{L_{2d}}$ which
  are invariant with respect to $\sigma_r$ with $r^2=-2$, because then
  $\sigma_r\in \Tilde{\Orth}^+(L_{2d})$. 

  (ii) The Fourier expansion of $\Phi_{12}$ at the standard
  0-dimensional cusp is defined by the hyperbolic unimodular lattice
  $L_{1,25}=U\oplus 3E_8(-1)$ (see (\ref{tube}) and (\ref{fourier})):
$$
\Phi_{12}(Z) =\sum_{u\in L_{1,25},\, (u,u)=0} \ a(u)\exp(2\pi i
(u,Z)). 
$$
The weight $12$ is singular, therefore the hyperbolic norm of the
index of any non-zero Fourier coefficient is zero. 

Let us fix a root $r\in R_l\subset L_{1,25}$ (any root is equivalent
to such a root).  We denote by $L_r$ the orthogonal complement of $r$
in $L_{1,25}$.  We have $Z=Z_r+ zr$, where $Z_r\in \cH (L_r)$ and
$z\in \CC$.  We note that $\Phi_{12}(Z_r)\equiv 0$. The function
$$
\Phi_r(Z_r)=\left.\frac{\Phi_{12}(Z)}{(Z, r)} \,\right\vert_{\cH(L_r)}
$$
is the first coefficient of the Taylor expansion of the function
$\Phi_{12}(Z_r+zr)$ in $z$. 

The summation in the Fourier expansion of $\Phi_r(Z_r)$ is taken over
the dual lattice $L_r^\vee$. We note that
$$
L_r \oplus \ZZ r\subset L_{1,25}\subset L_r^\vee 
\oplus \ZZ (r/2). 
$$ 
Let us calculate 
$$
\left.\frac{\partial \Phi_{12}(Z_r+zr)}{\partial z}\right\vert_{z=0}. 
$$
We get non-zero Fourier coefficient only for indices $u=u_r+m(r/2)$,
where $u_r\in L_r^\vee$ and $0\ne m\in \ZZ$.  In this case
$(u_r,u_r)=m^2/2>0$.  Thus the first derivative has non-zero Fourier
coefficient only for indices $u_r$ with positive square.  Doing this
for every $r$ we see that the Fourier expansion of $F_l$ at the
canonical cusp contains only indices with positive hyperbolic norm. 
Thus $F_l$ is a cusp form. 

The components of the branch divisor are divisors
$$
\cF_{2d}(r)=\pi(\{Z\in \cD_{L_{2d}}\mid (Z,r)=0\})
$$
defined by reflective vectors $r\in L_{2d}$, by
Corollary~\ref{reflK3}.  For a $(-2)$-vector $r\in L_{2d}$, the form
$F_l(Z)$ has a zero along $\cF_{2d}(r)$ (see (i)). 

Now we can finish the proof using Lemma~\ref{orthcomp2d}.  If $r$ is a
$(-2d)$-reflective vector and $L_r=r^\perp_{L_{2d}}$, then the divisor
$\cF_{2d}(r)$ coincides with the modular projection $\pi(\cD_{L_r})$
of the homogeneous domain of the lattice $L_r$ of signature $(2,18)$. 
According to Lemma~\ref{orthcomp2d}, $(L_r)^\perp_{L_{2,26}}$ is a
root lattice with $N\ge 112$ roots ($E_8$ has $240$ roots, $E_7$ has
$126$ and $D_8$ has $112$).  Therefore the Borcherds form $\Phi_{12}$
has a zero of order $N\ge 112>N_l$ along the subdomain $\cD_{L_r}$. 
Thus $F_l$ is zero along the corresponding divisor $\cF_{2d}(r)$. 
\end{proof}

According to Theorem~\ref{cuspform} and Theorem~\ref{general_gt} the
main point for us is the following.  We want to know for which $2d>0$
there exists a vector
\begin{equation}\label{orth2}
  l\in E_8,\ l^2=2d,\ l\ \text{ is orthogonal to at least $2$ and at most $12$ 
    roots.} 
\end{equation}
\begin{theorem}\label{mainineq} 
  Such a vector $l$ in $E_8$ does exist if one of two inequalities
\begin{equation}\label{mineq}
4N_{E_7}(2d)>28N_{E_6}(2d)+63N_{D_6}(2d)
\end{equation}
or
\begin{equation}\label{mineqd}
  5N_{E_7}(2d)>28N_{E_6}(2d)+63N_{D_6}(2d)+378N_{D_5}(2d)
\end{equation}
is valid, where $N_L(2d)$ denotes the number of representations of
$2d$ by the lattice $L$. 
\end{theorem}
\begin{proof}
  Let us fix a root $a \in E_8$. This choice gives us a realisation of
  the lattice $E_7$ as a sublattice of $E_8$:
$$
E_7\cong E_7^{(a)}=a^\perp_{E_8}. 
$$
We have the following decomposition of the set of roots $R(E_8)$:
$$
R(E_8)=R(E_7)\sqcup X_{114}\qquad\text{where }
X_{114}=\{c\in R(E_8)\mid c\cdot a \ne 0\}
$$
and $|X_{114}|=|R(E_8)|-|R(E_7)|=240-126=114$. 

\begin{lemma}\label{A2decomposition}
  The roots have the following properties:
\begin{itemize}
\item [(i)] $X_{114}$ is the union of $28$ root systems of type $A_2$
  such that $R(A_2^{(i)})\cap R(A_2^{(j)})=\{\pm a\}$ for any $i\ne
  j$. 
\item[(ii)] Let $A_2(a,c)\ne A_2(a,d)$ be two $A_2$-lattices generated
  by roots $a$, $c$ and $a$, $d$. Then
$$
A_3(a,c,d)=A_2(a,c)+A_2(a,d)
$$
is a lattice of type $A_3$ containing one and only one copy of $A_1$
from $E_7^{(a)}$. 
\item[(iii)] Let us take three different $A_2(a,c_i)$ ($i=1,2,3$). 
  Then their sum
$$
S=\sum_{i=1}^3 A_2(a,c_i)
$$ 
is a lattice of type $A_4$ or $D_4$.  The first one contains $20$
roots, the second contains $24$ roots.  In both cases exactly six
roots of $S$ are in $E_7^{(a)}$. 
\end{itemize}
\end{lemma}
\begin{proof} (i) Recall that $|b\cdot c|\le 2$ for any roots $b$,
  $c\in R(E_8)$.  If $b\cdot c=\pm 2$ then $b=\pm c$.  We can assume
  that $a\cdot c=-1$ (if not we replace $c$ by $-c$).  The lattice
  $A_2(a,c)=\ZZ a+\ZZ c$ is a lattice of $A_2$-type.  Any
  $A_2$-lattice contains six roots
$$
R(A_2(a,c))=\{\,\pm a,\ \pm c,\ \pm(a+c)\,\}. 
$$
$A_2(a,c)$ is generated by any pair of linearly independent roots. 
Therefore
$$
A_2(a,c_1)\cap A_2(a,c_2)=\{ \pm a\}
$$ 
if the root lattices are distinct.  \smallskip

(ii) $c\ne \pm d$ implies that $c\cdot d=0$ or $\pm 1$.  Suppose that
$c\cdot d=0$. Then the sum of the lattices is of type $A_3$ ($a\cdot
c= a\cdot d=-1$ and $c\cdot d=0$).  This lattice contains $12$ roots
$$
R(A_3(a,c,d))=\pm (a,\ c,\ d,\ a+c, \ a+d; \ a+c+d). 
$$ 
The first five roots are elements of $X_{114}$ and $a+c+d\in
E_7^{(a)}$. 

If $c\cdot d=1$ then $(a+d)\cdot c=0$ and we come back to the first
case.  If $c\cdot d=-1$ then $(a+d)\cdot c=-2$, $c=-(a+d)$ and
$A_2(a,c)= A_2(a,d)$.  \smallskip

(iii) As in the proof of $2)$ we can suppose that 
$c_1\cdot c_2=c_2\cdot c_3=0$ and $c_1\cdot c_3=0$ or $1$. 

If $c_1\cdot c_3=1$, then we see that  $S$ has a root basis of
type $A_4$. 
\vskip.5cm
\noindent\hskip-.5cm\begin{picture}(300,10)(55,10)
\put(120,10){\circle*{5}}
\put(115,0){$c_{3}$}
\put(120,10){\vector(1,0){62}}
\put(180,10){\circle*{5}}
\put(170,0){$-c_{1}$}
\put(180,10){\vector(1,0){62}}
\put(240,10){\circle*{5}}
\put(225,0){$a+c_{1}$}
\put(240,10){\vector(1,0){62}}
\put(300,10){\circle*{5}}
\put(295,0){$c_{2}$}
\end{picture}
\vskip1cm
$A_4$ has $20$ roots. They are
$$
\pm(a,\ c_i,\ a+c_i,\ a+c_1+c_2, \ a+c_2+c_3,\ c_1-c_3)
\quad\text{where }i=1,2,3. 
$$
Only the last three roots belong to $E_7^{(a)}$. 
\smallskip

If $c_1\cdot c_3=0$ then the roots $c_1$, $a$, $c_2$, $c_3$ form a
basis of $S$. In this case $S$ has type $D_4$ ($a\cdot c_i=-1$ for all
$i$ and the other scalar products are zero).  This root system
contains all roots of $A_4$ except $\pm(c_1-c_3)$ and the roots
$$
\pm(a+c_1+c_3,\ a+c_1+c_2+c_3,\ 2a+c_1+c_2+c_3). 
$$
The six roots from $E_7^{(a)}$ are $\pm(a+c_i+c_j)$. 
\end{proof}

Now we can finish the proof of Theorem~\ref{mainineq}.  Let us assume
that every $l\in E_7^{(a)}$ with $l^2=2d>0$ is orthogonal to at least
$14$ roots in $E_8$ including $\pm a$.  The others are some roots in
$E_7^{(a)}$ ($126$ roots), or in $X_{114}\setminus \{\pm a\}$ ($112$
roots).  If $l$ is orthogonal to $b\in X_{114}\setminus \{\pm a\} $
then $l$ is orthogonal to the lattice $A_2(a,b)$.  Therefore using
Lemma~\ref{A2decomposition} we have
\begin{equation}\label{union}
l\in \bigcup_{i=1}^{28} (A_2^{(i)})^\perp_{E_8} \cup
\bigcup_{j=1}^{63} (A_1^{(j)})^\perp_{E_7}. 
\end{equation}
We recall that $(A_2)^\perp_{E_8}\cong E_6$, $(A_1)^\perp_{E_7}\cong
D_6$ and $(A_1\oplus A_1)^\perp_{E_8}\cong D_6$.  Let denote by $n(l)$
the number of components in (\ref{union}) containing the vector $l$. 
We have calculated this vector exactly $n(l)$ times in the sum
$$
28N_{E_6}(2d)+63N_{D_6}(2d). 
$$
We shall consider several cases. 

(a). Suppose that $l\cdot c\ne 0$ for any $c\in X_{114}\setminus \{\pm
a\}$.  Then $l$ is orthogonal to at least $6$ copies of $A_1$ in
$E_7^{(a)}$ and $n(l)\ge 6$.  \smallskip

Now we suppose that there exist $c\in X_{114}\setminus \{\pm a\}$ such
that $l\cdot c=0$. Then $l$ is orthogonal to $A_2(a,c)$ which is one
of the $28$ subsystems of the bouquet $X_{114}$.  \smallskip

(b). If $l$ is orthogonal to only one $A_2^{(i)}$ ($6$ roots) then $l$
is orthogonal to at least $4$ copies of $A_1$ ($8$ roots) in
$E_7^{(a)}$.  Thus $n(l)\ge 5$.  \smallskip

(c). If $l$ is orthogonal to exactly two $A_2^{(i)}$ and $A_2^{(j)}$
in $X_{114}$ then $l$ is orthogonal to $A_3=A_2^{(i)}+A_2^{(j)}$
having $12$ roots and containing only one $A_1$ from $E_7^{(a)}$. Thus
$l$ is orthogonal to another $A_1$ in $E_7^{(a)}$. Therefore $n(l)\ge
4$.  \smallskip

(d). If $l$ is orthogonal to three or more $A_2^{(i)}$ then their sum
contains three $A_1\subset E_7^{(a)}$ and $n(l)\ge 6$.  \smallskip

We see that under our assumption $n(l)\ge 4$ for any $l\in E_7^{(a)}$. 
Therefore we have proved that if every $l\in E_7^{(a)}$ with $l^2=2d$
is orthogonal to at least $14$ roots then
$$
28N_{E_6}(2d)+63N_{D_6}(2d)\ge 4N_{E_7}(2d). 
$$
Moreover $n(l)$ can be equal to $4$ only in case~(c). In this case
$l\in (A_3)^\perp_{E_8}\cong D_5$ and there are $\binom{28}{2}=378$
pairs of $A_2$-subsystems in $X_{114}$. This gives us the second
inequality
$$
28N_{E_6}(2d)+63N_{D_6}(2d)\ge 5N_{E_7}(2d)-378N_{D_5}(2d). 
$$
\end{proof}

The inequalities (\ref{mineq}) and (\ref{mineqd}) fail only for a
finite number of $d$ because their left- and right-hand sides have the
asymptotics $O(d^{5/2})$ and $O(d^2)$. 

\begin{proposition}\label{P-ex} 
  A vector $l\in E_8$ satisfying the condition (\ref{orth2}) does
  exist if $d\not \in P_{ex}$, where
\begin{gather*}
  P_{ex}=\{\,1\le m\le 100\ (m\ne 96);\quad 101\le m\le 127\ (m \text{
    is odd});
  \\
  m=110,\ 131, 137,\ 143\,\}. 
\end{gather*}
\end{proposition}
\begin{proof}
  The Jacobi theta-series of the lattice $E_8$ coincides with the
  Jacobi-Eisenstein series $E_{4,1}(\tau,z)$ of weight $4$ and index
  $1$.  Let us fix a root $a\in E_8$. We have
$$
E_{4,1}(\tau,z)=\sum_{l\in E_8} \exp(\pi i\, l^2\tau+2\pi i\, l\cdot a z)=
1+\sum_{m\ge 1} e_{4,1}(m,n) \exp(2\pi m\tau+ nz). 
$$ 
$N_{E_7}(2m)=e_{4,1}(m,0)$, since the orthogonal complement of $a$ in
$E_8$ is $E_7$. 

The Fourier coefficients $e_{4,1}(m,n)$ were calculated in \cite{EZ}. 
In particular
$$
N_{E_7}(2m)=\frac{2^6\pi^3}{15}\frac{L^{Z}_{4m}(3)}{\zeta(3)}\,m^{5/2}
$$ 
where 
$$
L^{Z}_{D}(s)=\sum_{t\ge 1}\frac{\#\{\,x\mod 2t\mid x^2\equiv D\mod
  4t\,\}}{t^s}. 
$$ 
It is evident that $L^{Z}_{4m}(3)>9/8$ (one has to take only two terms
for $t=1$ and $t=2$). Thus
\begin{equation}\label{NE7estimate}
  N_{E_7}(2m)>\frac{24\pi^3}{5\zeta(3)}\,m^{5/2}>c(E_7) m^{5/2},
\end{equation}
where $c(E_7)=123.8$. In fact this estimate is quite good: a
computation with PARI shows that $N_{E_7}(314)\approx 124.73\times
(157)^{5/2}$

We can find simple exact formulae for $N_{E_6}(2m)$ and $N_{D_6}(2m)$. 
Let $\chi_3$ and $\chi_4$ be the unique non-trivial Dirichlet
characters modulo $3$ and $4$ respectively. For a Dirichlet character
$\chi$ we put
$$
\sigma_k(m,\chi)=\sum_{d|m} \chi(d)d^k,
\qquad 
\tilde\sigma_k(m,\chi)=\sum_{d|m} \chi\left(\frac{m}d\right)d^k. 
$$ 
\begin{lemma}\label{NE6andND6}
  The number of representations of $2m$ by the quadratic forms $E_6$
  and $D_6$ are
$$
\begin{aligned}
N_{E_6}(2m)&=81\tilde\sigma_2(m,\chi_3)-9\sigma_2(m,\chi_3),\\
N_{D_6}(2m)&=64\tilde\sigma_2(m,\chi_4)-4\sigma_2(m,\chi_4). 
\end{aligned}
$$
\end{lemma}
\begin{proof} 
  The second identity is well-known. This is the number of
  representations of $2m$ by six squares. To prove the first identity
  we consider the theta-series of $E_6$:
$$
\theta_{E_6}(\tau)=\sum_{l\in E_6}e^{\pi i (l\cdot l)}
\in M_3(\Gamma_0(3),\chi_3)=M_3(\Gamma_1(3)). 
$$
The dimension of $M_3(\Gamma_1(3))$ is equal to $2$.  We can construct
a basis with the help of Eisenstein series $G_k^{\alpha}$, where
$\alpha=(a,b)\in (\ZZ/N\ZZ)^2$,
$$
G_k^{\alpha}(\tau)=\sum_{(n,m)\equiv (a,b)\mod N}\ (n\tau+m)^{-k}. 
$$ 
Using the relation $G_k^{\alpha}|_k\gamma=G_k^{\alpha\gamma}$ (where
$\gamma\in \SL_2(\ZZ)$) for $k=3$ and $N=3$ we obtain two modular
forms in $M_3(\Gamma_1(3))$, namely $G_3^{(0,1)}$ and
$G_3^{(1,0)}+G_3^{(1,1)}+G_3^{(1,2)}$.  The Fourier expansion of
$G_k^{\alpha}$ was found by Hecke (see \cite{Kob}).  Normalising both
series we obtain a basis of $M_3(\Gamma_0(3),\chi_3)$ consisting of
$$ 
\begin{aligned} 
E_3^{(\infty)}(\tau,\chi_3)&=1-9\sum_{m\ge 1}\sigma_2(m,\chi_3)q^m,\\
E_3^{(0)}(\tau,\chi_3)&=\sum_{m\ge 1}\tilde\sigma_2(m,\chi_3)q^m
\qquad\qquad(q=e^{2\pi i \tau}). 
\end{aligned}
$$
We note that the first series is proportional to
$(\eta^3(\tau)/\eta(3\tau))^3$ and it vanishes at the cusp $0$. The
second series vanishes at $i\infty$.  The lattice $E_6$ has $72$
roots. Therefore
\begin{equation}\label{thE6}
\theta_{E_6}(\tau)=81E_3^{(0)}(\tau,\chi_3)+E_3^{(\infty)}(\tau,\chi_3). 
\end{equation}
This gives us the formula for $N_{E_6}(2m)$.  Applying the same method
to the theta-series $\theta_{D_6}\in M_3(\Gamma_0(4),\chi_4)$ we
obtain that
\begin{equation}\label{thD6}
\theta_{D_6}(\tau)=64E_3^{(0)}(\tau,\chi_4)+E_3^{(\infty)}(\tau,\chi_4),
\end{equation}
where
$$
\begin{aligned} 
E_3^{(\infty)}(\tau,\chi_4)&=1-4\sum_{m\ge 1}\sigma_2(m,\chi_4)q^m,\\
E_3^{(0)}(\tau,\chi_4)&=\sum_{m\ge 1}\tilde\sigma_2(m,\chi_4)q^m. 
\end{aligned}
$$
\end{proof}
Using these representations we can get an upper bound for
$N_{E_6}(2m)$ and $N_{D_6}(2m)$. It is clear that
$$
\sigma_2(m,\chi_3)
=\chi_3(m)\tilde\sigma_2(m,\chi_3)
\qquad\text{if }\ m\not \equiv 0\mod 3. 
$$
For any $C\equiv 1\mod 3$ we have the following bound
$$
\frac{\tilde\sigma_2(m,\chi_3)}{m^2}=\sum_{d|m}\frac{\chi_3(d)}{d^2}<
\sum_{1\le l\le C,\ l\equiv 1\mod 3} l^{-2}+
\bigg(\zeta(2)-\sum_{1\le n\le C+2,}n^{-2}\bigg). 
$$
Taking $C=19$ we get that for any $m$ not divisible by $3$
\begin{equation}\label{NE6estimate}
N_{E_6}(2m)=\tilde\sigma_2(m,\chi_3)(81-9\chi_3(m))<c(E_6)m^2,
\end{equation}
where $c(E_6)=103.69$. 

If $m=3^km_1$ then $\sigma_2(m,\chi_3)=\sigma_2(m_1,\chi_3)$, so
the last inequality is valid for any $m$.  For $D_6$ one can take
$C=21$ in a similar sum.  As a result we get
\begin{equation}\label{ND6estimate}
N_{D_6}(2m)<c(D_6)m^2,
\end{equation}
where $c(D_6)=75.13$. 

Using the estimates (\ref{NE7estimate}), (\ref{NE6estimate}) and
(\ref{ND6estimate}) for $N_L(2m)$, where $L=E_7$, $E_6$ and $D_6$, we
obtain that the main inequality (\ref{mineq}) of
Theorem~\ref{mainineq} is valid if
$$
m\ge 238>\left(\frac{28c(E_6)+63c(D_6)}{4c(E_7)}\right)^2. 
$$
For smaller $m$ we can use another formula for the theta-series of
$E_7$ (see \cite[(112)]{CS})
\begin{equation}\label{thetaE7}
\theta_{E_7}(\tau)
=\theta_3(2\tau)^7+7\theta_3(2\tau)^3\theta_2(2\tau)^4,
\end{equation}
where 
$$
\theta_3(2\tau)=\sum_{n=-\infty}^{\infty} q^{n^2},\qquad
\theta_2(2\tau)=\sum_{n=-\infty}^{\infty} q^{(n+\frac{1}2)^2}. 
$$
Moreover (see \cite[(87)]{CS})
\begin{equation}\label{thetaDn}
\theta_{D_n}(\tau)
=\frac{1}{2}(\theta_3(\tau)^n+\theta_3(\tau+1)^n). 
\end{equation}
Using (\ref{thetaE7}) and (\ref{thetaDn}) together with (\ref{thE6})
we can compute (using PARI) the first $240$ Fourier coefficients of
the function
$$
5\theta_{E_7}-28\theta_{E_6}-63\theta_{D_6}-378\theta_{D_5}. 
$$
The indices of the negative coefficients form the set $P_{ex}$ of
$d$ for which the inequality (\ref{mineqd}) of Theorem~\ref{mainineq}
fails. 
\end{proof}

Now we are going to analyse the main condition (\ref{orth2}) for some
$d\in P_{ex}$ from Proposition~\ref{P-ex}.  Moreover we are also
looking for vectors with $d\le 61$ orthogonal to exactly $14$ roots. 
Such vectors produce cusp forms $F_l$ of weight $19$ due to Theorem
\ref{cuspform}.  \smallskip

Let $e_i$ ($1\le i\le 8$) be a euclidean basis of the lattice $\ZZ^8$
($(e_i, e_j)=\delta_{ij}$).  We consider the Coxeter basis of simple
roots in $E_8$ (see \cite{Bou})

\noindent\hskip-.5cm
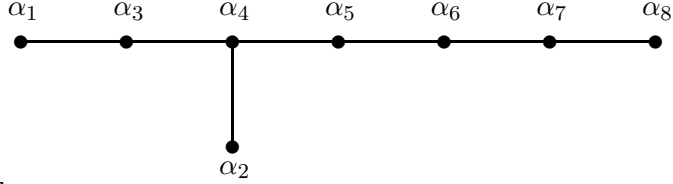
\begin{picture}(300,10)(55,10)
\put(100,0){\circle*{5}}
\put(95,10){$\alpha_1$}
\put(100,0){\vector(1,0){42}}
\put(140,0){\circle*{5}}
\put(135,10){$\alpha_3$}
\put(140,0){\vector(1,0){42}}
\put(180,0){\circle*{5}}
\put(175,10){$\alpha_4$}
\put(180,1){\vector(0,-1){43}}
\put(180,-40){\circle*{5}}
\put(175,-50){$\alpha_2$}
\put(180,0){\vector(1,0){42}}
\put(220,0){\circle*{5}}
\put(215,10){$\alpha_5$}
\put(220,0){\vector(1,0){42}}
\put(260,0){\circle*{5}}
\put(255,10){$\alpha_6$}
\put(260,0){\vector(1,0){42}}
\put(300,0){\circle*{5}}
\put(295,10){$\alpha_7$}
\put(300,0){\vector(1,0){42}}
\put(340,0){\circle*{5}}
\put(335,10){$\alpha_8$}
\end{picture}
\vskip2cm
where
\begin{gather*}
\alpha_1=\frac 1{2}(e_1+e_8)-\frac 1{2}(e_2+e_3+e_4+e_5+e_6+e_7),\\
\alpha_2=e_1+e_2,\quad \alpha_k=e_{k-1}-e_{k-2}\ \ (3\le k\le 8)
\end{gather*}
and $E_8=\langle\alpha_1,\dots \alpha_8\rangle_\ZZ$. 

Let $L_S=\langle \alpha_i\mid i\in S\rangle_\ZZ \subset E_8$ be a
sublattice of $E_8$ generated by some simple roots ($S\subset
\{1,\dots,8\}$).  We assume that $\#R (L_S)\le 12$, where $R(L_S)$ is
the set of roots of $L_S$.  We can find the orthogonal complement of
$L_S$ in $E_8$ using fundamental weights $\omega_j$, {i.e.} the basis
of $E_8$ dual to the basis $\{\alpha_i\}_{i=1}^8$.  We have
$$
L_S^\perp=(L_S)^\perp_{E_8}
=\langle \omega_j\mid  j\not\in S\rangle_\ZZ. 
$$
Any vector of $L_S^\perp$ is orthogonal to all roots of $L_S$.  If
$l\in L_S^\perp$ is orthogonal to an additional root $r$ of $E_8$
($r\not\in R(L_S)$) then we obtain a linear relation on the
coordinates of $l$ in the basis $\omega_j$ $(j\not\in S)$. 
Considering all roots of $E_8$ we can formulate a condition on the
coordinates of $l\in L^\perp_S$ to be orthogonal to at most $12$ roots
(or to exactly $14$ roots). We shall analyse four different lattices
$L_S$.  \medskip

\noindent{\bf I}. $L_1=4A_1$, $\ \# R(4A_1)=8$ and $L_1^\perp=4A_1$. 
\smallskip

We put
$$
L_1=\langle \alpha_2,\, \alpha_3,\, \alpha_5,\, \alpha_7\rangle_\ZZ 
=\langle e_2+e_1,\, e_2-e_1,\, e_4-e_3,\, e_6-e_5\rangle_\ZZ \cong 4A_1. 
$$ 
This root lattice $L_1$ gives us vectors of norm $2d$ for most $d\in
P_{ex}$.  $L_1$ is a primitive sublattice of $E_8$.  Therefore
$L_1^\perp$ is a lattice with the same discriminant form and
$L_1^\perp\cong 4A_1$. More exactly,
$$
L_1^\perp=\langle\,\omega_1,\, \omega_4,\, \omega_6,\, \omega_8\rangle_\ZZ
=\langle\, e_3+e_4,\, e_5+e_6,\, e_7+e_8,\, e_7-e_8\rangle_\ZZ. 
$$
This representation follows easily from the formulae for the
fundamental weights of $E_8$ (see \cite[Plat VII]{Bou}):
\begin{gather*}
\omega_2=\frac{1}{2}(e_1+\dots+e_7+5e_8),\quad
\omega_3=\frac{1}{2}(-e_1+e_2+\dots+e_7+7e_8),\\
\omega_k=e_{k-1}+\dots+e_7+(9-k)e_8\quad (4\le k\le 8), \quad \omega_1=2e_8. 
\end{gather*}
Any vector
\begin{equation}\label{linL1}
l=m_3(e_3+e_4)+m_5(e_5+e_6)+m_7(e_7+e_8)+m_8(e_7-e_8)\in L_1^\perp
\end{equation}
is orthogonal to $8$ roots of $L_1$.  The root system of $E_8$
contains $112$ integral and $128$ half-integral roots:
$$
\pm e_i\pm e_j\quad (i<j),\quad
\frac{1}{2}\sum_{i=1}^{8}(-1)^{\nu_i}e_i\quad
\text{with }\ \sum_{i=1}^{8} {\nu_i}\equiv 0\mod 2. 
$$
If $l$ is orthogonal to a half-integral root $r$ then 
\begin{multline}\label{lr1}
2(l\cdot r)=m_7((-1)^{\nu_7}+(-1)^{\nu_8})+m_8((-1)^{\nu_7}-(-1)^{\nu_8})+\\
m_3((-1)^{\nu_3}+(-1)^{\nu_4})+m_5((-1)^{\nu_5}+(-1)^{\nu_6})=0. 
\end{multline}
We note that only one of $m_7$ or $m_8$ appears.  Let us assume that
this identity contains three non-zero terms: $m_{7,8}\pm m_3\pm m_5
=0$ (by $m_{7,8}$ we mean $m_7$ or $m_8$).  Then $l$ is orthogonal to
$4$ additional half-integral roots.  There are two choices for
$(\nu_1,\nu_2)$ and one can change the sign of the root.  A similar
result, {i.e.} a relation $m_7\pm m_8\pm m_{3,5}=0$ and $4$ additional
integral roots, is obtained if $l$ is orthogonal to the integral roots
$e_{7,8}\pm e_{3,4}$ or $e_{7,8}\pm e_{5,6}$. 

If (\ref{lr1}) contains only two non-zero terms then we have a
relation of type $m_{7,8}\pm m_{3,5}=0$. In this case $l$ is
orthogonal to $8$ additional half-integral roots: there are two
choices for $(\nu_3,\nu_4)$ (or $(\nu_5,\nu_6)$), for $(\nu_1,\nu_2)$
and the change of the sign.  We can also have $m_{7,8}=0$, and then
the number of half-integral roots orthogonal to $l$ is equal to $16$. 

If $l$ is orthogonal to an integral root $r\not\in L_1$, which has not
been considered above, then we get a relation $m_3=\pm m_5$ or
$m_7=\pm m_8$ with $8$ additional roots or $m_{3,5}=0$ with $16$
additional integral roots.  For example, if $m_7=m_8$ then $l$ is
orthogonal to $\pm(e_8\pm e_{1,2})$; if $m_{3}=0$ then $l$ is
orthogonal to $\pm (e_{3,4}\pm e_{1,2})$.  Therefore we have proved
the following

\begin{proposition}\label{4A1} 
  $l\in L_1^\perp$ (see (\ref{linL1})) is orthogonal to at least $8$
  and at most $12$ roots of $E_8$ if and only if
\begin{itemize}
\item[(i)] $m_j\ne 0$ for any $j$ and $m_i\ne m_j$ for any $i\ne j$;
\item[(ii)] There is at most one relation of type $m_k=\pm m_i\pm m_j$
  for $i<j<k$. 
\end{itemize}
\end{proposition}
This lemma gives us a set of vectors $l\in L_1^\perp$ with
$$
l^2=2(m_3^2+m_5^2+m_7^2+m_8^2)=2d\in P_{ex}
$$
such that $l$ is orthogonal to $8$ or to $12$ roots of $E_8$.  We list
these vectors in table~{\bf I-8,12}. 

\begin{center}
\begin{tabular}{|c|c||c|c||c|c|}
\hline
\multicolumn{6}{|c|}
{{\bf I-8,12.}\quad $L_1=4A_1$,\quad $l=(m_3,m_5,m_7,m_8)\in L_1^\perp$}\\
\hline
$d$&$l$&$d$&$l$&$d$&$l$\\
\hline
$46$&$(1,2,4,5)$&$84$&$(1,3,5,7)$&$110$&$(1,3,6,8)$\\
\hline
$50$&$(1,2,3,6)$&$85$&$(1,2,4,8)$&$111$&$(1,2,5,9)$\\
\hline
$54$&$(2,3,4,5)$&$86$&$(3,4,5,6)$&$113$&$(2,3,6,8)$\\
\hline
$57$&$(1,2,4,6)$&$90$&$(1,2,6,7)$&$117$&$(1,4,6,8)$\\
\hline
$62$&$(1,3,4,6)$&$91$&$(1,4,5,7)$&$119$&$(2,3,5,9)$\\
\hline
$63$&$(1,2,3,7)$&$93$&$(2,3,4,8)$&$121$&$(1,2,4,10)$\\
\hline
$65$&$(2,3,4,6)$&$94$&$(1,2,5,8)$&$123$&$(1,3,7,8)$\\
\hline
$66$&$(1,2,5,6)$&$95$&$(1,3,6,7)$&$125$&$(3,4,6,8)$\\
\hline
$70$&$(1,2,4,7)$&$98$&$(2,3,6,7)$&$127$&$(1,3,6,9)$\\
\hline
$71$&$(1,3,5,6)$&$99$&$(3,4,5,7)$&$131$&$(3,4,5,9)$\\
\hline
$74$&$(2,3,5,6)$&$102$&$(1,2,4,9)$&$137$&$(2,4,6,9)$\\
\hline
$78$&$(1,2,3,8)$&$105$&$(1,2,6,8)$&$143$&$(1,5,6,9)$\\
\hline
$79$&$(1,2,5,7)$&$107$&$(1,3,4,9)$&&\\
\cline{1-4}
$81$&$(2,4,5,6)$&$109$&$(2,4,5,8)$&&\\
\hline
\end{tabular}
\end{center}

\medskip
\noindent{\bf II.} $L_2=2A_1\oplus A_2$, $\# R(2A_1\oplus A_2)=10$. 
\smallskip

Our second example is the sublattice 
$$
L_2=\langle \alpha_2,\, \alpha_3,\, \alpha_5,\, \alpha_6\rangle_\ZZ=
\langle e_2+e_1,\, e_2-e_1,\, e_4-e_3,\, e_5-e_4\rangle_\ZZ 
\cong 2A_1\oplus A_2. 
$$
Then using the dual basis $\omega_j$ we obtain that 
\begin{multline}\label{linL2}
L_2^\perp =\langle \omega_1, \omega_4,\omega_7,\omega_8\rangle=
\langle e_3+e_4+e_5+e_6, e_6+e_7, e_7-e_8, e_7+e_8\rangle\\
=\left\{l=m_5(e_3+e_4+e_5)+\sum_{i=6}^{8} m_ie_i\mid
m_5+m_6+m_7+m_8\text{ is even}\right\}. 
\end{multline}
We note that $L_2^\perp$ is not a root lattice. 

The vector $l$ is orthogonal to a half-integral root $r$ if 
$$
2(l\cdot r)=m_5((-1)^{\nu_3}+(-1)^{\nu_4}+(-1)^{\nu_5})+m_6(-1)^{\nu_6}+
m_7(-1)^{\nu_7}+m_8(-1)^{\nu_8}\!=0. 
$$
There are two different cases:
\begin{itemize}
\item[---] if $3m_5=\pm m_6\pm m_7\pm m_8$ then there are $4$
  half-integral roots orthogonal to $l$, since there are two choices
  for $(\nu_1,\nu_2)$ and for the sign of $r$;
\item[---] if $m_5=\pm m_6\pm m_7\pm m_8$ then there are $12$
  half-integral roots orthogonal to $l$, since there are three choices
  for $(\nu_3,\nu_4, \nu_5)$. 
\end{itemize}
Let us find integral roots of $E_8$ (not in $L_2$) orthogonal to $l$: 
\begin{itemize}
\item[---] if $m_i=0$ ($i=6$, $7$ or $8$) then there are $8$ roots
  $\pm(e_{1,2}\pm e_i)$;
\item[---] if $m_5=0$ then there are $24$ roots $\pm(e_{1,2}\pm
  e_{3,4,5})$;
\item[---] if $m_i=\pm m_5$ ($i=6$, $7$ or $8$) then there are $6$
  roots $\pm(e_i\mp e_{3,4,5})$;
\item[---] if $m_i=\pm m_j$ ($6\le i<j\le 8$) then there are $2$ roots
  $\pm(e_i\mp e_j)$. 
\end{itemize}
Therefore we obtain 
\begin{proposition}\label{2A1} 
  $l\in L_2^\perp$ (see (\ref{linL2})) is orthogonal to exactly $10$
  roots of $E_8$ if and only if
\begin{itemize}
\item[(i)] $m_j\ne 0$ for any $j$ and $m_i \ne \pm m_j$ for any $i<j$;
\item[(ii)] $km_5\ne \pm m_6\pm m_7\pm m_8$, where $k=1$ or $3$. 
\end{itemize}
Moreover $l\in L_2^\perp$ is orthogonal to exactly $14$ roots of $E_8$
if (i) and (ii) for $k=1$ are valid and there is exactly one relation
of type $3m_5=\pm m_6\pm m_7\pm m_8$. 
\end{proposition} 

Some $l\in L_2^\perp$ orthogonal to $10$ roots in $E_8$ and having
norm $l^2=3m_5^2+m_6^2+m_7^2+m_8^2=2d\in P_{ex}$ are given in
table~{\bf II-10}. 

\begin{center}
\begin{tabular}{|c|c||c|c||c|c|}
\hline
\multicolumn{6}{|c|}
{{\bf II-10.}\quad $L_2=2A_1\oplus A_2$,\quad $l=(m_5;\,m_6,m_7,m_8)\in L_2^\perp$}\\
\hline
$d$&$l$&$d$&$l$&$d$&$l$\\
\hline
$58$&$(1;\,2,3,10)$&$75$&$(6;\,1,4,5)$&$89$&$(2;\,6,7,9)$\\
\hline
$60$&$(3;\,2,5,8)$&$80$&$(3;\,4,6,9)$&$97$&$(4;\,1,8,9)$\\
\hline
$64$&$(5;\,1,4,6)$&$82$&$(5;\,3,4,8)$&$100$&$(7;\,1,4,6)$\\
\hline
$67$&$(2;\,4,5,9)$&$83$&$(2;\,1,3,12)$&$101$&$(4;\,1,3,12)$\\
\hline
$72$&$(3;\,1,4,10)$&$87$&$(6;\,1,4,7)$&$103$&$(8;\,1,2,3)$\\
\hline
$73$&$(4;\,3,5,8)$&$88$&$(1;\,2,5,12)$&$115$&$(4;\,1,9,10)$\\
\hline
\end{tabular}
\end{center}

The vectors from the tables~{\bf I-8,12} and {\bf II-10} produce cusp
forms $F_l(Z)$ of weights $16$, $18$ (table~{\bf I-8,12}) or $17$
(table~{\bf II-10}) for all $d>61$ in the set $P_{ex}$ except $d=68$,
$69$, $77$, $92$. 

The vectors from $L_2^\perp$ with $l^2=2d$ and $d\le 61$ that are
orthogonal to exactly $14$ roots of $E_8$ are given in
table~{\bf{II-14}}. 

\begin{center}
\begin{tabular}{|c|c||c|c||c|c|}
\hline
\multicolumn{6}{|c|}
{{\bf II-14.}\quad $L_2=2A_1\oplus A_2$,\quad $l=(m_5;\,m_6,m_7,m_8)\in L_2^\perp$}\\
\hline
$d$&$l$&$d$&$l$&$d$&$l$\\
\hline
$40$&$(1;\,2,3,8)$&$48$&$(3;\,1,2,8)$&$55$&$(4;\,1,5,6)$\\
\hline
$43$&$(2;\,1,3,8)$&$52$&$(1;\,2,4,9)$&$61$&$(2;\,1,3,10)$\\
\hline
\end{tabular}
\end{center}
\medskip

\noindent{\bf III}. $L_3=A_3$, $\# R(A_3)=12$. 
\smallskip

The root lattice $A_3$ is maximal. Therefore any sublattice of type
$A_3$ in $E_8$ is primitive. Analysing the discriminant form of the
orthogonal complement of $A_3$ we obtain that it is isomorphic to
$D_5$.  We put
$$
L_3=\langle\,\alpha_2,\,\alpha_4,\,\alpha_3\rangle_\ZZ=
\langle\,e_2+e_1,\,e_3-e_2,\,e_2-e_1\rangle_\ZZ\cong A_3. 
$$
Then 
$$
L_3^\perp=\bigg\{l=\sum_{i=4}^{8}m_ie_i\mid\sum_{i=4}^{8}m_i\equiv
0\mod 2\bigg\}\cong D_5. 
$$
As above we obtain 

\begin{proposition}\label{A3} $l\in L_3^\perp$ is orthogonal 
to exactly $12$ roots of $E_8$ if and only if
\begin{itemize}
\item[(i)] $m_j\ne 0$ for any $j$; 
\item[(ii)] $m_i \ne \pm m_j$ for any $i<j$;
\item[(iii)] $\sum_{i=4}^{8} \pm m_i\ne 0$ for any choice of the signs. 
\end{itemize}
Moreover $l\in L_3^\perp$ is orthogonal to exactly 
$14$ roots of $E_8$ if (i) and (iii) are valid and there is only 
one relation of type $m_i=\pm m_j$ for $4\le i<j\le 8$. 
\end{proposition} 

See table~{\bf III} for several vectors $l\in L_3^\perp$ orthogonal to
$N_l$ roots ($N_l=12$ or $14$) in $E_8$ and having norm
$l^2=\sum_{i=4}^{8}m_i^2=2d$. 

\begin{center}
\begin{tabular}{|c|c|c||c||c|c|}
\hline
\multicolumn{6}{|c|}
{{\bf III.}\quad $L_3=A_3$,\quad $l=(m_4,m_5,m_6,m_7,m_8)\in L_3^\perp$}\\
\hline
$d$&$l$&$N_l$&$d$&$l$&$N_l$\\
\hline
$69$&$(2,3,5,6,8)$&$12$&$53$&$(1,4,4,3,8)$&$14$\\
\hline
$42$&$(1,3,3,4,7)$&$14$&$54$&$(1,3,3,5,8)$&$14$\\
\hline
$48$&$(1,1,2,3,9)$&$14$&$56$&$(1,1,5,6,7)$&$14$\\
\hline
$49$&$(2,2,4,5,7)$&$14$&$59$&$(1,2,2,3,10)$&$14$\\
\hline
$51$&$(1,6,6,2,5)$&$14$&$63$&$(3,4,4,6,7)$&$14$\\
\hline
\end{tabular}
\end{center}
\medskip

\noindent{\bf IV}. $L_4=A_1\oplus A_2$, $\# R(A_1\oplus A_2)=8$. 
\smallskip

For any sublattice $A_1\oplus A_2$ in $E_8$ we see that its orthogonal
complement is isomorphic to $A_5$, since $(A_2)^\perp_{E_8}=E_6$ and
$(A_1)^\perp_{E_6}=A_5$.  We put
$L_4=\langle\,\alpha_1,\,\alpha_2,\,\alpha_3\rangle_\ZZ \cong
A_1\oplus A_2 $.  Then
$$
L_4^\perp=\bigg\{l=\sum_{i=3}^{8} m_ie_i \mid m_8=\sum_{i=3}^{7}m_i \bigg\}. 
$$
If $l$ is orthogonal to a half-integral root distinct from $\alpha_1$,
$\alpha_1+\alpha_3\in L_4$ then we get a relation of the form
$$
m_{i_1}+\dots+m_{i_k}=0,\quad\text{ where }\quad 3\le i_1<\dots i_k\le 7,\quad
1\le k\le 5. 
$$ 
If any relation of this type is valid then $l$ is orthogonal to $4$
additional {\it half-integral\/} roots.  Considering the scalar
products with {\it integral\/} roots we see that
\begin{itemize}
\item[---] if $m_i=0$ ($3\le i\le 8$) then $l$ is orthogonal to $8$
  roots $\pm(e_{1,2}\pm e_i)$;
\item[---] if $m_i=\pm m_j$ ($3\le i<j\le 8$) then $l$ is orthogonal
  to $2$ roots $\pm(e_i\mp e_j)$. 
\end{itemize}

We list some cases of these results in table~{\bf IV}. 
\begin{center}
\begin{tabular}{|c|c|c||c|c|c|}
\hline
\multicolumn{6}{|c|}
{{\bf IV.}\quad $L_4=A_1\oplus A_2$,\quad $l=(m_3, m_4,m_5,m_6,m_7;\,m_8)\in L_4^\perp$}\\
\hline
$d$&$l$&$N_l$&$d$&$l$&$N_l$\\
\hline
$68$&$(1,3,4,5,-7;\,8)$&$12$&$92$&$(1,1,2,3,5;\,12)$&$10$\\
\hline
$77$&$(2,3,4,5,-8;\,6)$&$12$&$40$&$(1,1,2,3,-8;\,-1)$&$14$\\
\hline
\end{tabular}
\end{center}
\medskip

It is possible to formulate a result for this case analogous to
Propositions~\ref{4A1}, \ref{2A1} and \ref{A3}, but we do not need it. 

An extensive computer search for vectors $l$ orthogonal to at least~$2$ and
at most~$14$ roots for other $d\in P_{ex}$ has not found any. 

Now we have everything we need to prove our main theorem,
Theorem~\ref{mainthm}. For $d>61$ and for $d=46$, $50$, $54$, $57$,
$58$, $60$ there exists a vector $l$ satisfying
condition~(\ref{orth2}), either by Proposition~\ref{P-ex} or listed in
one of the tables. Hence Theorem~\ref{cuspform} provides us with a
suitable cusp form of low weight. Since the dimension of $\cF_{2d}$ is
$19$, Theorem~\ref{main_sings_theorem} guarantees the existence of a
compactification with only canonical singularities and hence
Theorem~\ref{mainthm} follows by using the low weight cusp form trick,
according to Theorem~\ref{general_gt}. 

If $d$ is not as above but $d\ge 40$ and $d\ne 41$, $44$, $45$, $47$
then we have a cusp form of weight~$19$ arising from a vector $l$
orthogonal to $14$ roots, listed in one of the tables. This gives rise
to a canonical form and hence, by Freitag's result, the Kodaira
dimension of $\cF_{2d}$ is non-negative, as stated in
Theorem~\ref{general_gt}.

\bibliographystyle{alpha}

\bigskip
\noindent
V.A.~Gritsenko\\
Universit\'e Lille 1\\
Laboratoire Paul Painlev\'e\\
F-59655 Villeneuve d'Ascq, Cedex\\
France\\
{\tt valery.gritsenko@math.univ-lille1.fr}
\bigskip

\noindent
K.~Hulek\\
Institut f\"ur Algebraische Geometrie\\
Leibniz Universit\"at Hannover\\
D-30060 Hannover\\ 
Germany\\
{\tt hulek@math.uni-hannover.de}
\bigskip

\noindent
G.K.~Sankaran\\
Department of Mathematical Sciences\\
University of Bath\\
Bath BA2 7AY\\
England\\
{\tt gks@maths.bath.ac.uk}

\end{document}